%% file: paper.tex
\documentclass{article}
\usepackage{jlcode}
\RequirePackage{amsmath}
\RequirePackage{amsfonts}
\RequirePackage{amssymb}
\RequirePackage{mathtools}  % fixes some problems in amsmath
\RequirePackage{mathrsfs}   % different caligraphic fonts via \mathscr{N}
\RequirePackage{siunitx}
\RequirePackage{bbm}        % \mathbbm{N}
\RequirePackage{nicefrac}   % compact symbols for 1/2, etc. (when it comes to units)
\RequirePackage[hyphens]{url}
\usepackage{arxiv}
\usepackage{kthcolors}
\usepackage{algorithm}
\usepackage{algpseudocode}

\algnewcommand\algorithmicinput{\hspace*{1.0em}\textbf{input:}}
\algnewcommand\Input{\item[\algorithmicinput]}
\algnewcommand\algorithmicoutput{\hspace*{1.0em}\textbf{output:}}
\algnewcommand\Output{\item[\algorithmicoutput]}

\newcommand{\algorithmicand}{\textbf{ and }}
\renewcommand{\And}{\algorithmicand}

\providecommand{\given}{}
\DeclarePairedDelimiter{\parentheses}{(}{)}
\DeclarePairedDelimiter{\brackets}{[}{]}
\DeclarePairedDelimiter{\braces}{\lbrace}{\rbrace}

\DeclarePairedDelimiter{\Verts}{\lVert}{\rVert}

\providecommand{\SetSymbol}[1][]{%
  \nonscript\:#1\vert
  \allowbreak
  \nonscript\:
  \mathopen{}%
}
\DeclarePairedDelimiterX{\Set}[1]{\lbrace}{\rbrace}{
  \renewcommand{\given}{\SetSymbol[\delimsize]}
  #1
}

\providecommand{\expect}[2][{}k]{\ensuremath{\operatorname{\mathbb{E}}_{#1}\brackets*{#2}}}
\providecommand{\var}[2][{}k]{\ensuremath{\operatorname{Var}}_{#1}\brackets*{#2}}

\providecommand{\norm}[1]{\ensuremath{\Verts*{#1}}}
\providecommand{\argmin}[2][x]{\ensuremath{\arg\,\min_{#1}\braces*{#2}}}

\providecommand{\sbj}{\ensuremath\text{subject to}}
\providecommand{\st}{\ensuremath\,\text{s.t.}\,}

\DeclareMathOperator*{\minimize}{minimize}

\usepackage{pgfplots}
\usepackage{tikz}

\pgfplotsset{compat=1.14}

\usetikzlibrary{arrows.meta,intersections,fit,math,shapes,calc,positioning}

\begin{document}

\title{Efficient Stochastic Programming in Julia}

\author{Martin Biel \\
  Division of Decision and Control Systems\\
  School of EECS, KTH Royal Institute of Technology\\
  SE-100 44 Stockholm, Sweden\\
  \texttt{mbiel@kth.se} \\
  \And Mikael Johansson \\
  Division of Decision and Control Systems\\
  School of EECS, KTH Royal Institute of Technology\\
  SE-100 44 Stockholm, Sweden\\
  \texttt{mikaelj@kth.se}}

\maketitle

\begin{abstract}
  \sloppy
  We present \jlinl{StochasticPrograms.jl}, a user-friendly and powerful open-source framework for stochastic programming written in the Julia language. The framework includes both modeling tools and structure-exploiting optimization algorithms. Stochastic programming models can be efficiently formulated using expressive syntax and models can be instantiated, inspected, and analyzed interactively. The framework scales seamlessly to distributed environments. Small instances of a model can be run locally to ensure correctness, while larger instances are automatically distributed in a memory-efficient way onto supercomputers or clouds and solved using parallel optimization algorithms. These structure-exploiting solvers are based on variations of the classical L-shaped and progressive-hedging algorithms. We provide a concise mathematical background for the various tools and constructs available in the framework, along with code listings exemplifying their usage. Both software innovations related to the implementation of the framework and algorithmic innovations related to the structured solvers are highlighted. We conclude by demonstrating strong scaling properties of the distributed algorithms on numerical benchmarks in a multi-node setup.
\end{abstract}

\section{Introduction}
\label{sec:introduction}

Stochastic programming is an effective mathematical framework for modeling multi-stage decision problems that involve uncertainty~\cite{Birge2011}. It has been used to model complex real-world problems in diverse fields such as power systems~\cite{Fleten2007,Groewe-Kuska2005,petra_real-time_2014}, finance~\cite{Krokhmal2005,Zenios2005}, and transportation~\cite{Powell1987,Powell2005}. The classical setting is linear stochastic programs where an actor takes decisions in two stages:
\begin{equation*}
  \text{initial decision } x \quad\rightarrow\quad \text{observation } \omega \quad\rightarrow\quad \text{recourse action } y(x, \xi(\omega))
\end{equation*}
The actor first takes a decision $x$. Then, the realization of a random event $\omega$ alters the state of the world. The actor can observe $\omega$ and take a recourse action $y$ with respect to $x$ and the output of some random variable $\xi(\omega)$. We are interested in finding the optimal decision $x$, accounting for the ability to make a recourse action once $\omega$ has been observed. The notion of an optimal decision is captured by letting $x$ and $y$ be optimization variables in linear programs, where $\xi(\omega)$ parameterizes the  second-stage problem for each event $\omega$.

In applications, a stochastic program models some real-world decision problem under a statistical model of $\xi$. We can then compute approximations of optimal decision policies by solving approximated instances of the stochastic program. In brief, this involves computing a first-stage decision $\hat{x}$ that is optimal in expectation over a set of second-stage scenarios $\xi(\omega_i)$ sampled from the model of $\xi$. This technique is known as sampled average approximation (SAA). In the linear setting, one can in principle formulate sampled instances on a extensive form that considers all available scenarios at once. This mathematical program can be solved using standard linear programming solvers, including both open-source solvers such as GLPK~\cite{glpk} and commercial solvers such as Gurobi~\cite{gurobi}. However, the size of the extensive form grows linearly in the number of scenarios, and industry-scale applications typically involve 10,000+ scenarios. For example, the 24-hour unit commitment problem studied in~\cite{petra_real-time_2014} has 16,384 scenarios and the resulting extensive form has 4 billion variables. Solving the extensive form in such applications becomes practically infeasible. Moreover, the memory requirement for storing the stochastic program instances will eventually exceed the capacity of a single machine. This clarifies the need for a distributed approach for modeling large-scale stochastic programs. Structure-exploiting decomposition methods~\cite{Rockafellar1991,van_slyke_l-shaped_1969} that operate in parallel on distributed data become essential to solve large-scale instances.

\subsection{Contribution}
\label{sec:contribution}

\sloppy
In this work, we present a user-friendly open-source software framework for efficiently modeling and solving stochastic programs in a distributed-memory setting. The framework allows researchers to formulate complex stochastic models and quickly typeset and test novel optimization algorithms. Stochastic programming educators will benefit from the clean and expressive syntax. Also, the framework supports analysis tools and stochastic programming constructs from classical theory and leading textbooks. Industrial practitioners can make use of the framework to rapidly formulate complex models, analyze small instances locally, and then run large-scale instances in production on a supercomputer or a cloud cluster. We implemented the framework in the Julia~\cite{Bezanson2017} programming language. Henceforth, we refer to the framework as SPjl. The framework is freely available through the registered Julia package \jlinl{StochasticPrograms.jl}.

The design philosophy adopted during implementation of SPjl is centered around flexibility and efficiency, with the aim to provide a feature-rich and user-friendly experience. Also, the framework should be scalable to support large-scale problems. With this in mind, we adhered to the fundamental principle that the optimization modeling should be separated from the data modeling. This design principle results in two key software innovations: deferred model instantiation and data injection. Optimization models are formulated in stages using a straightforward syntax that simultaneously specifies the data dependencies between the stages. The data structures related to future scenarios, and their statistical properties, are defined separately. An essential consequence of this design is that we can efficiently distribute stochastic program instances in memory, reducing interprocess communication to a minimum. Many computations involving distributed stochastic programs can then natively be run in parallel. Moreover, when the sample space is infinite, it becomes possible to adequately distinguish between the abstract representation of a stochastic program and finite sampled instances. The design also enables swift implementation of various constructs from classical stochastic programming theory. Another design choice is that the solver suites included in the framework are developed using policy-based techniques. We have shown in prior work how policy-based design can be used to create customizable and efficient optimization algorithms~\cite{polojl}. In short, SPjl is a powerful, versatile, and extensible framework for stochastic programming. It provides both an educational setting for new practitioners and a research setting were experts can further the field of stochastic programming.

\sloppy
We developed SPjl in Julia, which has several distinct benefits. Through just-in-time compilation and type inference, Julia can achieve C-like performance while being as expressive as similar dynamic languages such as Python or Matlab. Using the high-level metaprogramming capabilities of Julia, it is possible to create domain-specific tools with expressive syntax and high performance. Another benefit is access to Julia's large and rapidly expanding ecosystem of libraries, many of which play a central role in SPjl. For example, the parallel capabilities of SPjl are implemented using the standard library module for distributed computing, while optimization models are formulated using the JuMP~\cite{Dunning2017} ecosystem. JuMP is an algebraic modeling language implemented in Julia using similar metaprogramming tools. It has been shown to achieve similar performance to AMPL~\cite{Dunning2017}, with syntax that is both readable and expressive. Also, it is possible to mutate model instances at runtime, which we utilize in the structure-exploiting algorithms. Recently, the backend of JuMP was redesigned into the new MathOptInterface~\cite{moi}. The redesign introduces automatic reformulation bridges, which are used frequently in the current implementation of the SPjl framework. JuMP implements interfaces to many third-party optimization solvers, both open-source and commercial. These can be hooked in to solve extensive forms of stochastic programs or subproblems that arise in decomposition methods.

\subsection{Related work}
\label{sec:related-work}

We give a short survey of similar software packages and highlight distinguishing features of SPjl. The most similar approach is the PySP framework~\cite{Watson2012}, implemented in the Python language. Optimization models in PySP are created using Pyomo~\cite{Hart2017}; an algebraic modeling language also implemented in Python. In contrast, SPjl is written in the Julia language and formulates optimization models in JuMP, which has been shown to outperform Pyomo in various benchmarks~\cite{Dunning2017}. In PySP, stochastic programs are composed of multiple \texttt{.dat} files and \texttt{.py} files, and the models are solved by running different solver scripts. In SPjl, all models are described in pure Julia and can be created, analyzed and solved in a single interactive session. Moreover, all operations are natively distributed in memory and run in parallel if multiple Julia processes are available. The parallel capabilities of PySP extend to running parallelized versions of the solver scripts. The primary function of PySP is to formulate and solve stochastic programs, while SPjl also provides a large set of stochastic programming constructs and analysis tools. The expressiveness of the modeling syntax can be compared by observing how the well-known farmer problem~\cite{Birge2011} is modeled using PySP~\cite{Watson2012} and how it is modeled using SPjl, as shown in Appendix~\ref{sec:farmer-problem}. In particular, the PySP definition requires about $\num{100}$ lines of code spread out over four different files, while SPjl requires $\num{30}$ lines of code with the added benefit of being more readable. In addition, the resulting model can be analyzed interactively in Julia in a user-friendly way.

A more extensive list of similar software approaches is provided in~\cite{Watson2012}, along with comparisons to PySP. This allows for a transitive comparison to SPjl. Other notable examples include the commercial FortSP solvers~\cite{fortsp} coupled with the AMPL extension SAMPL for modeling. Out of all these approaches, SPjl has the most user-friendly interface and is also freely available.

The StructJuMP package~\cite{Huchette2014} provides a simple interface to create block-structured JuMP models. The primary reason for developing StructJuMP was to facilitate a parallel modeling interface to existing structured solvers~\cite{lubin_parallel_2013, petra_real-time_2014} that operate in computer clusters. These parallel solvers are implemented in C++ and are parallelized using MPI. This led to StructJuMP also making use of MPI to distribute stochastic programs in blocks. Apart from formulating distributed stochastic programs in a cluster, StructJuMP does not offer any modeling tools nor any way to generate the extensive form of a stochastic program. In comparison, SPjl provides numerous analysis tools as well as a compatible suite of structured solvers. In addition, SPjl natively distributes and solves stochastic programs using Julia, without relying on external software such as MPI.

\section{Preliminaries}
\label{sec:preliminaries}

We give a short mathematical introduction to linear stochastic programming. The purpose is to provide background for the code examples presented in the subsequent section and also to keep this work self-contained. A more thorough introduction to the field is given in the textbook by~\cite{Birge2011}.

\subsection{Stochastic programming}
\label{sec:stoch-progr-1}

A linear two-stage recourse model enables a simple but powerful framework for making decisions under uncertainty. We formalize this procedure in the following brief review. The first-stage decision made by the actor is denoted by $x \in \mathbb{R}^n$. We associate $x$ with a linear cost function $c^Tx$ that the actor pays after making the decision. Moreover, $x$ is constrained to the standard polyhedron in linear programming, i.e.
\begin{equation*}
  \Set{x \in \mathbb{R}^n \given Ax = b,\; x \geq 0}
\end{equation*}
where $A \in \mathbb{R}^{p \times n}$ and $b \in \mathbb{R}^p$. The recourse actions are represented by $y \in \mathbb{R}^m$. To describe the uncertainty in the decision problem, we consider some probability space $(\Omega,\mathcal{F},\pi)$ where $\Omega$ is a sample space, $\mathcal{F}$ is a $\sigma$-algebra over $\Omega$ and $\pi: \mathcal{F} \to [0,1]$ is a probability measure. Let $\xi(\omega): \Omega \to \mathbb{R}^{N}$ be some random variable on $\Omega$ and let $\operatorname{\mathbb{E}}_{\xi}$ denote expectation with respect to $\xi$. We can now let $\omega \in \Omega$ denote a scenario observed after deciding $x$. The scenario affects both cost and the constraints of the recourse action. Specifically, after realization of $\omega$, the following second-stage problem is formulated to determine $y$ with respect to $x$ and $\xi(\omega)$:
\begin{equation} \label{eq:lssubprob}
  \begin{aligned}
    Q(x,\xi(\omega)) = \min_{\mathclap{y \in \mathbb{R}^m}} & \quad q_{\omega}^T y \\
    \st & \quad T_{\omega}x + Wy = h_{\omega} \\
    & \quad y \geq 0.
  \end{aligned}
\end{equation}
In other words, the random variable takes on the form $\xi(\omega) = \begin{pmatrix}
  q_{\omega} & T_{\omega} & h_{\omega}
\end{pmatrix}^T$ in this linear setting. Note that $q_{\omega} \in \mathbb{R}^m$, $T_{\omega} \in \mathbb{R}^{q \times n}$ and $h_{\omega} \in \mathbb{R}^q$ are scenario-dependent while $W \in \mathbb{R}^{q \times m}$ is fixed. This is a standard setting in literature, which covers a wide range of problems~\cite{Birge2011}. It is  possible to define $W$ as scenario-dependent in the framework, but standard algorithms are then no longer certain to  converge. Now, we formulate the two-stage recourse problem as follows.
\begin{equation} \label{eq:linearsp}
  \begin{aligned}
    \minimize_{\mathclap{x \in \mathbb{R}^n}} & \quad c^T x + \expect[\xi]{Q(x,\xi(\omega))} \\
    \sbj & \quad Ax = b \\
    & \quad x \geq 0,
  \end{aligned}
\end{equation}
The optimal value of~\eqref{eq:linearsp} is referred to as the \emph{value of the recourse problem} (VRP).

Apart from solving~\eqref{eq:linearsp}, we can compute two classical measures of stochastic performance. The first measures the value of knowing the random outcome before making the decision. This is achieved by taking the expectation in~\eqref{eq:linearsp} outside the minimization, to obtain the wait-and-see problem:
\begin{equation} \label{eq:waitsee}
  \mathrm{EWS} = \expect[\xi]{
    \begin{aligned}
      \min_{x \in \mathbb{R}^n} & \quad c^T x + Q(x,\xi(\omega)) \\
      \st & \quad Ax = b \\
      & \quad x \geq 0.
    \end{aligned}}
\end{equation}
Now, the first- and second-stage decisions are taken with knowledge about the uncertainty. The difference between the expected wait-and-see value and the value of the recourse problem is known as the \emph{expected value of perfect information}:
\begin{equation} \label{eq:evpi}
    \mathrm{EVPI} = \mathrm{EWS} - \mathrm{VRP}.
\end{equation}
The EVPI measures the expected loss of not knowing the exact outcome beforehand. It quantifies the value of having access to an accurate forecast.

Finally, we introduce the concept of decision evaluation to quantify the performance of a candidate first-stage decision $x$ in the stochastic program~\eqref{eq:linearsp}. The \emph{expected result} of $x$ is given by
\begin{equation} \label{eq:expectedresult}
  V(x) = c^Tx + \expect[\xi]{Q(x,\xi(\omega))}.
\end{equation}
This concept is used to compute the second classical measure. If the expectation in~\eqref{eq:linearsp} is instead taken inside the second-stage objective function $Q$, we obtain the expected-value-problem:
\begin{equation} \label{eq:ev}
  \begin{aligned}
    \minimize_{x \in \mathbb{R}^n} & \quad c^T x + Q(x,\expect[\xi]{\xi(\omega)}) \\
    \sbj & \quad Ax = b \\
    & \quad x \geq 0.
  \end{aligned}
\end{equation}
The solution to the expected-value-problem is known as the \emph{expected value decision}, and is denote by $\bar{x}$. The \emph{expected result} of taking the \emph{expected value decision} is known as the \emph{expected result of the expected value decision}:
\begin{equation} \label{eq:eev}
  \mathrm{EEV} = c^T \bar{x} + \expect[\xi]{Q(\bar{x},\xi(\omega))}.
\end{equation}
The difference between the value of the recourse problem and the expected result of the expected value decision is known as the \emph{value of the stochastic solution}:
\begin{equation} \label{eq:vss}
  \mathrm{VSS} = \mathrm{EEV} - \mathrm{VRP}.
\end{equation}
The VSS measures the expected loss of ignoring the uncertainty in the problem. A large VSS indicates that the second stage is sensitive to the stochastic data.

The EVPI, VSS, and VRP are important tools when gauging the performance of a stochastic model. All of these introduced measures are readily computed in the SPjl framework, which allows for easy analysis of user-defined models. Next, we discuss how to calculate the VRP, EVPI, and VSS depending on the form of the sample space $\Omega$.

\subsection{The finite extensive form and sample average approximation}
\label{sec:finite-extens-form-saa}

If $\Omega$ is finite, say with $n$ scenarios of probability $\pi_1, \dots, \pi_n$, then we can represent~\eqref{eq:linearsp} compactly as
\begin{equation} \label{eq:finitesp}
  \begin{aligned}
    \minimize_{\mathclap{x \in \mathbb{R}^n, y_s \in \mathbb{R}^m}} & \quad c^T x + \sum_{s = 1}^{n} \pi_s q_s^T y_s & \\
    \sbj & \quad Ax = b & \\
    & \quad T_s x + W y_s = h_s, \quad &&s = 1,\dots,n \\
    & \quad x \geq 0, \, y_s \geq 0, \quad &&s = 1,\dots,n.
  \end{aligned}
\end{equation}
We refer to this problem as the \emph{finite extensive form}. It is often recognized in literature as the \emph{deterministic equivalent problem} (DEP). Similar closed forms can be determined for the EVPI and the VSS. For small $n$, it is viable to solve this problem with standard linear programming solvers. For large $n$, decomposition approaches are required. In SPjl, the user provides a description of the abstract stochastic model~\eqref{eq:linearsp} and a separate description of the uncertainty model of $\xi$. These are then combined internally to generate instances of the finite form~\eqref{eq:finitesp}, which are stored and solved efficiently on a computer or a compute cluster.

If $\Omega$ is not finite, the stochastic program~\eqref{eq:linearsp} is exactly computable only under certain assumptions~\cite{Birge2011}. However, it is possible to formulate computationally tractable approximations of~\eqref{eq:linearsp} using the finite form~\eqref{eq:finitesp}. The most common approximation technique is the \textit{sample average approximation} (SAA)~\cite{saa}. Assume that we sample $n$ scenarios $\omega_s,\; s = 1,\dots,n$ independently from $\Omega$ with equal probability. These scenarios now constitute a finite sample space $\tilde{\Omega}$ and we can use them to create a sampled model in finite extensive form~\eqref{eq:finitesp}. An optimal solution to this sampled model approximates the optimal solution to~\eqref{eq:linearsp} in the sense that the empirical average second-stage cost $V_n = \frac{1}{n}\sum_{s = 1}^{n}q_s^T \hat{y}_s,$ where $\hat{y}_s = \argmin[y \in \mathbb{R}^m]{Q(x, \xi(\omega_s))}$, converges pointwise with probability $1$ to $\hat{V} = \expect[\xi]{Q(x,\xi(\omega))}$ as $n$ goes to infinity~\cite{saaconvergence}. Further, under certain assumptions it can be shown that $\sqrt{n}(V_n - \hat{V}) \to N(0,\var[\xi]{Q(\hat{x},\xi)})$ in distribution as $n$ goes to infinity~\cite{saadist}. This result provides a basis for calculating confidence intervals around the VRP of~\eqref{eq:linearsp}~\cite{saa, saacomp}, as well as around the EVPI and the VSS.

\subsection{Structure-exploiting solvers}
\label{sec:struct-expl-solv}

Efficient methods for storing and solving finite stochastic programs on the form~\eqref{eq:finitesp} are key for high-performance stochastic programming computations. Therefore, this has been a main focus in the development of the SPjl framework. An important insight is that the finite extensive form~\eqref{eq:finitesp} lends itself to block-decomposition approaches, which allow the stochastic program to be efficiently distributed in memory. Moreover, structure-exploiting solvers can be employed to solve the decomposed models efficiently. These approaches also readily extend to parallel settings where the stochastic program is distributed over several compute nodes. A key idea in the SPjl framework is to let the storage of the stochastic program depend on the type of optimizer used to solve it. In this way, the memory structure is optimized for the solver operation, and there is no redundant storage for other operations such as decision evaluation. We say that the underlying structure of the stochastic program is induced by the solver. Henceforth, we will refer to the treatment of~\eqref{eq:finitesp} as one large optimization problem as the \emph{deterministic} structure. This is the default structure for standard third-party solvers. For block-decomposition approaches, we adopt the terminology introduced in~\cite{Watson2012} and divide such strategies into two classes. In short, ``\emph{vertical} strategies decompose a stochastic program by stages'' while ``\emph{horizontal} strategies decompose a stochastic program by scenarios''~\cite{Watson2012}. In the following, we will introduce two different solver algorithms that fall into these two categories and highlight the stochastic program structures they induce.

\subsubsection{The L-shaped algorithm}
\label{sec:l-shaped-algorithm}

The L-shaped algorithm is an efficient cutting-plane method for solving the finite extensive form~\eqref{eq:finitesp} by decomposing into a master problem and a set of subproblems. The master problem has the form
\begin{equation*}
  \begin{aligned}
    \minimize_{\mathclap{x \in \mathbb{R}^n, y_s \in \mathbb{R}^m}} & \quad c^T x + \theta \\
    \sbj & \quad Ax = b & \\
    & \quad \theta \geq \tilde{Q}(x) \\
    & \quad x \geq 0,
  \end{aligned}
\end{equation*}
where $\tilde{Q}(x)$ is a lower bound on
\begin{equation*}
  Q(x) = \sum_{s = 1}^{n}\pi_sQ_s(x).
\end{equation*}
Here, each $Q_s(x)$ is the optimal value to a subproblem of the form~\eqref{eq:lssubprob}. The idea of the L-shaped algorithm is to generate increasingly tight piecewise linear lower bounds on $Q$. We refer to the memory structure inferred by the L-shaped algorithm henceforth as the \emph{vertical} structure.

During the L-shaped procedure, solution candidates $x_k$ are generated by solving the master problem~\eqref{eq:lsmaster}, which are then used to parameterize subproblems of the form~\eqref{eq:lssubprob}. Optimal dual variables in these subproblems are then used to improve the bound of $Q(x)$ before the next solution candidate $x_{k+1}$ is computed. Specifically, it follows from duality theory that $\lambda_s^T(h_s-T_sx)$, where $\lambda_s$ is the dual optimizer of~\eqref{eq:lssubprob}, is a valid support function for $Q_s(x)$, and hence,
\begin{equation*}
  \sum_{s = 1}^{n}\pi_s \lambda_s^T(h_s-T_sx)
\end{equation*}
is a valid support function for $Q(x)$. In the original formulation of the L-shaped algorithm~\cite{van_slyke_l-shaped_1969}, the above result is used at each iteration $k$ to construct \emph{optimality cuts} by introducing
\begin{equation*}
  \partial Q_{k} = \sum_{s = 1}^{n} \pi_s\lambda_s^T T_s \qquad q_k = \sum_{s = 1}^{n}\pi_s\lambda_s^Th_s,
\end{equation*}
and add to the master problem as the constraint $\partial Q_k x + \theta \geq q_k$. Aggregating the results from all subproblems in this way is known as the single-cut approach. This was later extended to a multi-cut variant where separate cuts are constructed for each subproblem~\cite{Birge1988}. If the iterate $x_k$ is not second-stage feasible, some subproblems will be infeasible. We handle this by solving the auxilliary problem:
\begin{equation} \label{eq:feascheck}
  \begin{aligned}
    \minimize_{\mathclap{y_s \in \mathrm{R}^m}} & \quad w_s = e^Tv_s^{+} + e^Tv_s^{-} \\
    \sbj & \quad Wy_s  + v_s^{+} - v_s^{-}= h_s - T_s x_k \\
    & \quad y_s \geq 0, \; v_s^{+} \geq 0, \; v_s^{-} \geq 0.
  \end{aligned}
\end{equation}
If $w_s > 0$, then subproblem $s$ is infeasible for the current iterate $x_k$. Further, it follows from duality theory that $\sigma_s^T(h_s-T_sx) \leq 0$, where $\sigma_s$ is the dual optimizer of~\eqref{eq:feascheck}, is necessary for $x$ to be second-stage feasible. The above result can be used to both check for second-stage infeasibility and construct \emph{feasibility cuts} by introducing
\begin{equation*}
  F_k = \begin{pmatrix}
    \sigma_1^T T_1 \\
    \vdots \\
    \sigma_f^T T_f
  \end{pmatrix}, \quad
  f_k = \begin{pmatrix}
    \sigma_1^Th_1 \\
    \vdots \\
    \sigma_f^T h_f
  \end{pmatrix}
\end{equation*}
for all infeasible subproblems $1,\dots,f$. Because $W$ has a finite number of bases, finitely many feasibility cuts are required to completely describe the set of feasible first-stage decisions~\cite{van_slyke_l-shaped_1969}. The optimality cuts and the feasibility cuts enter the master problem as follows:
\begin{equation} \label{eq:lsmaster}
  \begin{aligned}
    \minimize_{\mathclap{x \in \mathbb{R}^n}} & \quad c^T x + \theta \\
    \sbj & \quad Ax = b \\
    & \quad F_k x \geq f_k, \quad &&  \forall k \\
    & \quad \partial Q_k x + \theta \geq q_k, \quad && \forall k \\
    & \quad x \geq 0.
  \end{aligned}
\end{equation}
The master problem is then re-solved to generate the next iterate $x_{k+1}, \theta_{k+1}$. This is repeated until the gap between the upper bound $Q(x_k)$ and lower bound $\theta_{k+1}$ becomes small, upon which the algorithm terminates. Many variations can be introduced to improve the performance of the L-shaped algorithm. We provide an overview of such improvements available in SPjl in Section~\ref{sec:algor-impr}.

\subsubsection{The progressive-hedging algorithm}
\label{sec:progr-hedg-algor}

The progressive-hedging algorithm was first introduced in~\cite{Rockafellar1991}. In contrast to the L-shaped algorithm, applying progressive-hedging to solve~\eqref{eq:finitesp} yields a complete decomposition over the $n$ scenarios. The method is a specialization of the proximal-point algorithm~\cite{Rockafellar1976}, and convergence in the linear case~\eqref{eq:finitesp} is derived in~\cite{Rockafellar1991}. The main idea behind this approach is to introduce individual first-stage decisions $x_s$ to each scenario but force them to be equal. We then relax (dualize) these consistency constraints and solve the corresponding augmented Lagrangian problem. In other words, we consider the following problem:
\begin{equation} \label{eq:phform}
\begin{aligned}
  \minimize_{\mathclap{x_s \in \mathbb{R}^n, y_s \in \mathbb{R}^m}} & \quad \sum_{s = 1}^{n} \pi_s \parentheses*{c^Tx_s + q_s^T y_s} & \\
  \sbj & \quad x_s = \xi \quad &&s = 1,\dots,n \\
 & \quad Ax_s = b \quad &&s = 1,\dots,n \\
 & \quad T_s x_s + W y_s = h_s, \quad &&s = 1,\dots,n \\
 & \quad x_s \geq 0, \; y_s \geq 0, \quad &&s = 1,\dots,n.
\end{aligned}
\end{equation}
The consistency constraints $x_s = \xi,\; s = 1,\dots,n$ are called \emph{non-anticipative} because they make the $x_s$ independent of scenario and enforce the fact that the first-stage decision is known when the second-stage uncertainty is realized. We refer to the memory structure inferred by the progressive-hedging algorithm henceforth as the \emph{horizontal} structure. Separability across the $n$ scenarios is achieved by introducing the following regularized relaxation of each subproblem:
\begin{equation*}
\begin{aligned}
  \minimize_{\mathclap{x_s \in \mathbb{R}^n, y_s \in \mathbb{R}^m}} & \quad c^Tx_s + q_s^T y_s + \rho_s(x_s-\xi) + \frac{r}{2}\norm{x_s-\xi}_2^2 \\
  \sbj & \quad Ax_s = b \\
 & \quad T_s x_s + W y_s = h_s \\
 & \quad x_s \geq 0, \; y_s \geq 0.
\end{aligned}
\end{equation*}
The algorithm now proceeds by iteratively alternating between generating new admissible solutions $x_s^k, \; s = 1,\dots, n$, and an implementable solution $\xi_k$. In the two-stage setting, an admissible solution is feasible in every scenario, and an implementable solution is consistent in the sense that $x_s = \xi$ for all $s$. We obtain the implementable solution through aggregation:
\begin{equation*}
  \xi_k = \sum_{s = 1}^{n}\pi_sx_s^k
\end{equation*}
and the Lagrange multipliers are updated scenario-wise through
\begin{equation*}
  \rho_s^{k+1} = \rho_s^{k} + r(x_s^k-\xi_k).
\end{equation*}
Hence, the non-anticipative constraints are enforced while the dual variables converge. Progressive-hedging is a primal dual algorithm that is run until both the primal gap $\norm{\xi_k - \xi_{k-1}}_2^2$ and the dual gap $\sum_{s = 1}^{n}\pi_s\norm{x_s^k - \xi_k}_2^2$ are small.

\section{StochasticPrograms.jl}
\label{sec:spjl}

In this section, we showcase the capabilities of SPjl. We first give a brief overview of the framework and introduce the main functionality through a set of simple examples. Accompanying code excerpts are included. Next, we exemplify the effectiveness of SPjl model creation by giving a compact definition of the farmer problem. Finally, we summarize the algorithmic improvements and variations included in the framework.

\sloppy
SPjl extends the well-known JuMP syntax to support the definition of stages, decision variables, and uncertain parameters. Models are defined using the \jlinl{@stochastic_model} macro. This creates a lightweight model object that can be used to instantiate finite stochastic programs by supplying a description of the uncertain parameters. Specifically, the user provides a list of discrete scenarios, or a sampler object capable of generating scenarios, to the model object. The object then combines the model definition with the supplied uncertainty data and generates a finite stochastic program instance. The instantiated stochastic program can then be inspected, analyzed and solved in an interactive Julia session. This is useful in educational settings, but also for reasoning about complex models on a small scale. SPjl also supports reading problems specified in the SMPS format. For large-scale instances, SPjl provides scalable block-structured instantiation and structure-exploiting solvers that can operate in parallel. In addition, operations such as EWS calculation and decision evaluation are embarrassingly parallel over the subproblems. In other words, the workload is readily decoupled into independent subtasks that can be executed in parallel. This is leveraged when instantiating vertical or horizontal structures in distributed environments.

\sloppy
SPjl can be installed directly from the command line through Julia's package manager (\jlinl{pkg> add StochasticPrograms}). Provided that a basic linear quadratic solver, such as \jlinl{GLPK} or \jlinl{Ipopt}, is installed, all code examples in this paper can be repeated by copying the lines verbatim. A more extensive introduction to the framework is given by the ``Quick start'' section of the online documentation~\footnote{\url{https://martinbiel.github.io/StochasticPrograms.jl/dev/}}.

\subsection{A simple textbook example}
\label{sec:simple-example}

Consider the following simple instance of~\eqref{eq:linearsp}
\begin{equation} \label{eq:simpleproblem}
\begin{aligned}
 \minimize_{\mathclap{x_1,x_2 \in \mathbb{R}}} & \quad 100x_1 + 150x_2 + \expect[\omega]{Q(x_1,x_2,\xi(\omega))} \\
 \sbj & \quad x_1+x_2 \leq 120 \\
 & \quad x_1 \geq 40 \\
 & \quad x_2 \geq 20
\end{aligned}
\end{equation}
where
\begin{equation} \label{eq:simplesubproblem}
\begin{aligned}
 Q(x_1,x_2,\xi(\omega)) = \max_{\mathclap{y_1,y_2 \in \mathbb{R}}} & \quad q_1(\omega)y_1 + q_2(\omega)y_2 \\
 \st & \quad 6y_1+10y_2 \leq 60x_1 \\
 & \quad 8y_1 + 5y_2 \leq 80x_2 \\
 & \quad 0 \leq y_1 \leq d_1(\omega) \\
 & \quad 0 \leq y_2 \leq d_2(\omega)
\end{aligned}
\end{equation}
and the stochastic variable
\begin{equation*}
  \xi(\omega) = \begin{pmatrix}
  q_1(\omega) & q_2(\omega) & d_1(\omega) & d_2(\omega)
  \end{pmatrix}^T
\end{equation*}
parameterizes the second-stage model. This is a recurring textbook example and correctness of our numerical results can be verified by comparing with~\cite{Birge2011}.

In SPjl, we create the model~\eqref{eq:simpleproblem} in two steps. First, we formulate the optimization models as shown in Listing~\ref{lst:simpledef}.
\begin{lstlisting}[language = julia, float, floatplacement = H, caption={Definition of~\eqref{eq:simpleproblem} in SPjl.}, label={lst:simpledef}]
# Load SPjl framework
julia> using StochasticPrograms
# Create simple stochastic model
julia> simple_model = @stochastic_model begin
  @stage 1 begin
      @decision(model, x₁ >= 40)
      @decision(model, x₂ >= 20)
      @objective(model, Min, 100*x₁ + 150*x₂)
      @constraint(model, x₁+x₂ <= 120)
  end
  @stage 2 begin
      @uncertain q₁ q₂ d₁ d₂
      @recourse(model, 0 <= y₁ <= d₁)
      @recourse(model, 0 <= y₂ <= d₂)
      @objective(model, Max, q₁*y₁ + q₂*y₂)
      @constraint(model, 6*y₁ + 10*y₂ <= 60*x₁)
      @constraint(model, 8*y₁ + 5*y₂ <= 80*x₂)
  end
end;
\end{lstlisting}
\noindent
\sloppy
This creates a stochastic model where the two stages are given by the mathematical programs~\eqref{eq:simpleproblem} and~\eqref{eq:simplesubproblem}, expressed using an enhanced JuMP syntax. The \jlinl{@decision} and \jlinl{@recourse} lines work as standard \jlinl{@variable} definitions in JuMP, but behind the scenes they also specify internal data dependencies between the first and second stage; and the \jlinl{@uncertain} line annotates the random parameters and defines a point of data injection. The code specifies how the optimization models should be defined, but the actual model instantiation is deferred until we add a stochastic model of the uncertainties. We will consider two different distributions of $\xi$ and use the same model object \jlinl{simple_model} from Listing~\ref{lst:simpledef} to instantiate stochastic programs. This is a key feature in SPjl. The underlying stochastic model~\eqref{eq:linearsp} object can be re-used to generate different finite stochastic program instances. Regardless of the distribution of $\xi$, a stochastic program instance is always a finite program of the form~\eqref{eq:finitesp}. This allows us to evaluate the same problem under different uncertainty models and to automatically adapt the underlying memory structure to optimize solver performance.

\subsection{Finite sample space}
\label{sec:finite-sample-space}

First, let $\xi$ be a discrete distribution, taking on the values
\begin{equation*}
  \xi_1 = \begin{pmatrix}
    500 & 100 & 24 & 28
  \end{pmatrix}^T, \quad \xi_2 = \begin{pmatrix}
    300 & 300 & 28 & 32
  \end{pmatrix}^T
\end{equation*}
with probability $0.4$ and $0.6$ respectively. In Listing~\ref{lst:simplediscrete}, an instance of the stochastic program~\eqref{eq:simpleproblem} is created for this distribution.
\begin{lstlisting}[language = julia, float, caption = {Instantiation of~\eqref{eq:simpleproblem}.}, label = {lst:simplediscrete}]
# Create two scenarios
julia> ξ₁ = @scenario q₁ = 24.0 q₂ = 28.0 d₁ = 500.0 d₂ = 100.0 probability = 0.4;
       ξ₂ = @scenario q₁ = 28.0 q₂ = 32.0 d₁ = 300.0 d₂ = 300.0 probability = 0.6;
# Instantiate without optimizer
julia> sp = instantiate(simple_model, [ξ₁, ξ₂])
Stochastic program with:
 * 2 decision variables
 * 2 scenarios of type Scenario
Structure: Deterministic equivalent
Solver name: No optimizer attached.
# Print to show structure of generated problem
julia> print(sp)
Deterministic equivalent problem
Min 100 x₁ + 150 x₂ - 9.6 y₁₁ - 11.2 y₂₁ - 16.8 y₁₂ - 19.2 y₂₂
Subject to
 y₁₁ ≥ 0.0
 y₂₁ ≥ 0.0
 y₁₂ ≥ 0.0
 y₂₂ ≥ 0.0
 y₁₁ ≤ 500.0
 y₂₁ ≤ 100.0
 y₁₂ ≤ 300.0
 y₂₂ ≤ 300.0
 x₁ ∈ Decisions
 x₂ ∈ Decisions
 x₁ ≥ 40.0
 x₂ ≥ 20.0
 x₁ + x₂ ≤ 120.0
 -60 x₁ + 6 y₁₁ + 10 y₂₁ ≤ 0.0
 -80 x₂ + 8 y₁₁ + 5 y₂₁ ≤ 0.0
 -60 x₁ + 6 y₁₂ + 10 y₂₂ ≤ 0.0
 -80 x₂ + 8 y₁₂ + 5 y₂₂ ≤ 0.0
Solver name: No optimizer attached.
\end{lstlisting}
\noindent
\sloppy
This code uses the model recipe created in Listing~\ref{lst:simpledef} to create second-stage models for each of the supplied scenarios. Here, we have used the default scenario constructor \jlinl{@scenario}, where data values are named in accordance with the \jlinl{@uncertain} annotation. The deterministic structure (extensive form) is used by default. Because this is a small example, correctness of the generated problem is easily verified. We can now set an optimizer and solve the model, as shown in Listing~\ref{lst:depsolve}.
\begin{lstlisting}[language = julia, float, caption = {Solving the finite extensive form of~\eqref{eq:simpleproblem}.}, label = {lst:depsolve}]
julia> using GLPK
# Set the optimizer to GLPK
julia> set_optimizer(sp, GLPK.Optimizer)
# Optimize (deterministic structure)
julia> optimize!(sp)
# Check termination status
julia> @show termination_status(sp);
termination_status(sp) = MathOptInterface.OPTIMAL
# Query optimal value
julia> @show objective_value(sp);
objective_value(sp) = -855.833333333333
# Calculate EVPI
julia> EVPI(sp)
662.916666666667
# Calculate VSS
julia> VSS(simple_model, SimpleSampler(μ, Σ))
286.9166666666688
\end{lstlisting}
\noindent
\sloppy
The underlying memory structure can be set explicitly by setting the \jlinl{instantiation} keyword to any of the supported structures during model instantiation. Alternatively, if an optimizer is chosen during instantiation, an appropriate structure is chosen automatically. For example, if we instantiate the same problem with an L-shaped optimizer the vertical structure is used instead, as can be seen in Listing~\ref{lst:simplels}.
\begin{lstlisting}[language = julia, float, caption = {Re-instantiation and optimization of~\eqref{eq:simpleproblem} with an L-shaped optimizer}, label = {lst:simplels}]
# Instantiate with L-shaped optimizer
julia> sp = instantiate(simple_model, [ξ₁, ξ₂], optimizer = LShaped.Optimizer)
Stochastic program with:
 * 2 decision variables
 * 2 scenarios of type Scenario
Structure: Vertical
Solver name: L-shaped with disaggregate cuts
# Print to compare structure of generated problem
julia> print(sp)
First-stage
==============
Min 100 x₁ + 150 x₂
Subject to
 x₁ ∈ Decisions
 x₂ ∈ Decisions
 x₁ ≥ 40.0
 x₂ ≥ 20.0
 x₁ + x₂ ≤ 120.0

Second-stage
==============
Subproblem 1 (p = 0.40):
Max 24 y₁ + 28 y₂
Subject to
 y₁ ≥ 0.0
 y₂ ≥ 0.0
 y₁ ≤ 500.0
 y₂ ≤ 100.0
 x₁ ∈ Known
 x₂ ∈ Known
 6 y₁ + 10 y₂ - 60 x₁ ≤ 0.0
 8 y₁ + 5 y₂ - 80 x₂ ≤ 0.0

Subproblem 2 (p = 0.60):
Max 28 y₁ + 32 y₂
Subject to
 y₁ ≥ 0.0
 y₂ ≥ 0.0
 y₁ ≤ 300.0
 y₂ ≤ 300.0
 x₁ ∈ Known
 x₂ ∈ Known
 6 y₁ + 10 y₂ - 60 x₁ ≤ 0.0
 8 y₁ + 5 y₂ - 80 x₂ ≤ 0.0
Solver name: L-shaped with disaggregate cuts
# Set GLPK optimizer for the solving master problem and subproblems
julia> set_optimizer_attribute(sp, MasterOptimizer(), GLPK.Optimizer)
julia> set_optimizer_attribute(sp, SubproblemOptimizer(), GLPK.Optimizer)
# Optimize (vertical structure)
julia> optimize!(sp)
L-Shaped Gap  Time: 0:00:02 (6 iterations)
  Objective:       -855.8333333333358
  Gap:             0.0
  Number of cuts:  8
  Iterations:      6
# Check termination status and query optimal value
julia> @show termination_status(sp);
termination_status(sp) = MathOptInterface.OPTIMAL
julia> @show objective_value(sp);
objective_value(sp) = -855.8333333333358
\end{lstlisting}
\noindent
\sloppy
The same stochastic program has now been decomposed into a first-stage master problem and two second-stage subproblems. For completeness we also exemplify how the same problem is instantiated and solved using the progressive-hedging algorithm in Listing~\ref{lst:simpleph}.
\begin{lstlisting}[language = julia, float, caption = {Re-instantiation and optimization of~\eqref{eq:simpleproblem} with a progressive-hedging optimizer}, label = {lst:simpleph}]
# Instantiate with progressive-hedging optimizer
julia> sp = instantiate(simple_model, [ξ₁, ξ₂],
                        optimizer = ProgressiveHedging.Optimizer)
Stochastic program with:
 * 2 decision variables
 * 2 scenarios of type Scenario
Structure: Horizontal
Solver name: Progressive-hedging with fixed penalty
# Print to compare structure of generated problem
julia> print(sp)
Horizontal scenario problems
==============
Subproblem 1 (p = 0.40):
Min 100 x₁ + 150 x₂ - 24 y₁ - 28 y₂
Subject to
 y₁ ≥ 0.0
 y₂ ≥ 0.0
 y₁ ≤ 500.0
 y₂ ≤ 100.0
 x₁ ∈ Decisions
 x₂ ∈ Decisions
 x₁ ≥ 40.0
 x₂ ≥ 20.0
 x₁ + x₂ ≤ 120.0
 -60 x₁ + 6 y₁ + 10 y₂ ≤ 0.0
 -80 x₂ + 8 y₁ + 5 y₂ ≤ 0.0

Subproblem 2 (p = 0.60):
Min 100 x₁ + 150 x₂ - 28 y₁ - 32 y₂
Subject to
 y₁ ≥ 0.0
 y₂ ≥ 0.0
 y₁ ≤ 300.0
 y₂ ≤ 300.0
 x₁ ∈ Decisions
 x₂ ∈ Decisions
 x₁ ≥ 40.0
 x₂ ≥ 20.0
 x₁ + x₂ ≤ 120.0
 -60 x₁ + 6 y₁ + 10 y₂ ≤ 0.0
 -80 x₂ + 8 y₁ + 5 y₂ ≤ 0.0
Solver name: Progressive-hedging with fixed penalty
julia> using Ipopt
# Set Ipopt optimizer for soving emerging subproblems
julia> set_optimizer_attribute(sp, SubproblemOptimizer(), Ipopt.Optimizer)
# Silence Ipopt
julia> set_optimizer_attribute(sp, RawSubproblemOptimizerParameter("print_level"), 0)
# Optimize (horizontal structure)
julia> optimize!(sp)
Progressive Hedging Time: 0:00:05 (303 iterations)
  Objective:   -855.5842547490254
  Primal gap:  7.2622997706326046e-6
  Dual gap:    8.749063651111478e-6
  Iterations:  302
# Check termination status and query optimal value
julia> @show termination_status(sp);
termination_status(sp) = MathOptInterface.OPTIMAL
julia> @show objective_value(sp);
objective_value(sp) = -855.5842547490254
\end{lstlisting}

\subsection{Infinite sample space}
\label{sec:infin-sample-space}

To demonstrate how SPjl handles continuous distributions for uncertain parameters, we assume that the uncertainties in our simple example follow a multivariate normal distribution, $\xi \sim \mathcal{N}(\mu, \Sigma)$. In general, there is no closed form solution of~\eqref{eq:linearsp} when $\xi$ has a continuous distribution. However, by the law of large numbers, a viable discrete approximation can be obtained by sampling scenarios from the continuous distribution. In SPjl, we achieve this by creating a sampler object associated with the defined scenario structure. In Listing~\ref{lst:simplecont}, a sampler object for a multivariate distribution with
\begin{equation*}
  \mu = \begin{pmatrix}
 24 \\
 32 \\
 400 \\
 200
\end{pmatrix}, \quad \Sigma = \begin{pmatrix}
 2 & 0.5 & 0 & 0 \\
 0.5 & 1 & 0 & 0 \\
 0 & 0 & 50 & 20 \\
 0 & 0 & 20 & 30
\end{pmatrix}
\end{equation*}
is created and used to generate an instance of~\eqref{eq:simpleproblem} with $100$ sampled scenarios.
\begin{lstlisting}[language = julia, float, caption = {Creating a sampled instance of~\eqref{eq:simpleproblem} in SPjl.}, label = {lst:simplecont}]
julia> using Distributions
# Define sampler object
julia> @sampler SimpleSampler = begin
    N::MvNormal # Normal distribution

    SimpleSampler(μ, Σ) = new(MvNormal(μ, Σ))

    @sample Scenario begin
        # Sample from normal distribution
        x = rand(sampler.N)
        # Create scenario matching @uncertain annotation
        return @scenario q₁ = x[1] q₂ = x[2] d₁ = x[3] d₂ = x[4]
    end
end
# Create mean
julia> μ = [24, 32, 400, 200];
# Create variance
julia> Σ = [2 0.5 0 0
            0.5 1 0 0
            0 0 50 20
            0 0 20 30];
# Instantiate sampled stochastic program with 100 scenarios
julia> sp = instantiate(simple_model, SimpleSampler(μ, Σ), 100)
Stochastic program with:
 * 2 decision variables
 * 100 scenarios of type Scenario
Structure: Deterministic equivalent
Solver name: No optimizer attached.
\end{lstlisting}
\noindent
\sloppy
Note that the same stochastic model object defined in Listing~\ref{lst:simpledef} is used in Listing~\ref{lst:simplecont} to generate the sampled instance.

With the ability to instantiate sampled models with an arbitrary number of scenarios, we can adopt the SAA methodologies developed in~\cite{saa} to calculate confidence intervals around the optimal value of~\eqref{eq:simpleproblem} as well as around the EVPI and the VSS. This is exemplified in Listing~\ref{lst:saasolve}.
\begin{lstlisting}[language = julia, float, caption = {Approximately solving~\eqref{eq:simpleproblem} when $\xi$ follows a normal distribution.}, label = {lst:saasolve}]
# Set optimizer to SAA
julia> set_optimizer(simple_model, SAA.Optimizer)
# Emerging stochastic programming instances solved by GLPK
julia> set_optimizer_attribute(simple_model, InstanceOptimizer(), GLPK.Optimizer)
# Set attributes that value solution speed over accuracy
julia> set_optimizer_attribute(simple_model, NumEvalSamples(), 300)
# Set target relative tolerance of the resulting confidence interval
julia> set_optimizer_attribute(simple_model, RelativeTolerance(), 5e-2)
# Approximate optimization using sample average approximation
julia> optimize!(simple_model, SimpleSampler(μ, Σ))
SAA gap Time: 0:00:03 (4 iterations)
  Confidence interval:  Confidence interval (p = 95%): [-1095.65 − -1072.36]
  Relative error:       0.021487453807842415
  Sample size:          64
# Check termination status
julia> @show termination_status(simple_model);
termination_status(sp) = MathOptInterface.OPTIMAL
# Query optimal value
julia> @show objective_value(simple_model);
objective_value(simple_model) = Confidence interval (p = 95%): [-1095.65 − -1072.36]
# Disable logging
julia> set_optimizer_attribute(simple_model, MOI.Silent(), true)
# Calculate approximate EVPI
julia> EVPI(simple_model, SimpleSampler(μ, Σ))
Confidence interval (p = 99%): [32.96 − 144.51]
# Calculate approximate VSS
julia> VSS(simple_model, SimpleSampler(μ, Σ))
Warning: VSS is not statistically significant to the chosen confidence level and tolerance
Confidence interval (p = 95%): [-0.05 − 0.05]
\end{lstlisting}
\noindent
\sloppy
These methods require re-solving sampled stochastic programs multiple times and the accuracy of the solution is increased by increasing the number of scenarios in the sampled models. Consequently, the parallel capabilities of SPjl become significant as these subproblems can become too large for single-core approaches. If multiple Julia processes are available, either locally or remotely, then the code in Listing~\ref{lst:simplecont} would automatically distribute the stochastic program on the available nodes in either a vertical or a horizontal structure. Although not practically required for this small example, this leads to significant performance gains for large-scale industrial models. See for example the scaling results presented in Section~\ref{sec:numerical-benchmarks}.

\section{The farmer problem}
\label{sec:farmer-problem}

To exemplify functional correctness, and allow for comparisons with similar tools, we consider the instructive farmer problem by~\cite{Birge2011}. Listing~\ref{lst:farmer} shows a suggested code excerpt for how the farmer problem can be defined in SPjl and Listing~\ref{lst:farmer_solve} shows how the problem can be instantiated, solved, and analyzed using various solvers. The correctness of the numerical values can be verified in~\cite{Birge2011}. For comparison, the same problem is defined in PySP as outlined in~\cite{Watson2012} in about $\num{100}$ lines spread out in separate files. Again, we stress that only $\num{30}$ lines of Julia code are required to define the farmer problem in SPjl. Moreover, the optimal value, as well as the EVPI and VSS, can be calculated interactively in the same Julia session. This feature distinguished SPjl from other similar tools such as PySP. The time required to solve the farmer problem using the L-shaped algorithm was $\num{0.57}$ seconds for both SPjl and PySP, measured on the same master node as the numerical benchmarks presented in the paper. Hence, there is no performance decrease from using SPjl instead of PySP for this small problem, with the added benefit of SPjl being more user-friendly.
\begin{lstlisting}[language = julia, float, caption={Definition of the farmer problem in SPjl}, label={lst:farmer}]
farmer = @stochastic_model begin
    @stage 1 begin
        @parameters begin
            Crops = [:wheat, :corn, :beets]
            Cost = Dict(:wheat=>150, :corn=>230, :beets=>260)
            Budget = 500
        end
        @decision(model, x[c in Crops] >= 0)
        @objective(model, Min, sum(Cost[c]*x[c] for c in Crops))
        @constraint(model, sum(x[c] for c in Crops) <= Budget)
    end
    @stage 2 begin
        @parameters begin
            Crops = [:wheat, :corn, :beets]
            Required = Dict(:wheat=>200, :corn=>240, :beets=>0)
            PurchasePrice = Dict(:wheat=>238, :corn=>210)
            SellPrice = Dict(:wheat=>170, :corn=>150, :beets=>36, :extra_beets=>10)
        end
        @uncertain ξ[c in Crops]
        @recourse(model, y[p in setdiff(Crops, [:beets])] >= 0)
        @recourse(model, w[s in Crops ∪ [:extra_beets]] >= 0)
        @objective(model, Min, sum(PurchasePrice[p] * y[p] for p in setdiff(Crops, [:beets]))
                   - sum(SellPrice[s] * w[s] for s in Crops ∪ [:extra_beets]))
        @constraint(model, minimum_requirement[p in setdiff(Crops, [:beets])],
            ξ[p] * x[p] + y[p] - w[p] >= Required[p])
        @constraint(model, minimum_requirement_beets,
            ξ[:beets] * x[:beets] - w[:beets] - w[:extra_beets] >= Required[:beets])
        @constraint(model, beets_quota, w[:beets] <= 6000)
    end
end
\end{lstlisting}

\begin{lstlisting}[language = julia, float, caption={Instantiation, optimization, and analysis of the farmer problem in SPjl}, label={lst:farmer_solve}]
# Define the three yield scenarios
julia> Crops = [:wheat, :corn, :beets];
       ξ₁ = @scenario ξ[c in Crops] = [3.0, 3.6, 24.0] probability = 1/3;
       ξ₂ = @scenario ξ[c in Crops] = [2.5, 3.0, 20.0] probability = 1/3;
       ξ₃ = @scenario ξ[c in Crops] = [2.0, 2.4, 16.0] probability = 1/3;
# Instantiate with GLPK optimizer
julia> farmer_problem = instantiate(farmer_model, [ξ₁,ξ₂,ξ₃], optimizer = GLPK.Optimizer)
# Optimize stochastic program (through extensive form)
julia> optimize!(farmer_problem)
# Inspect optimal decision
julia> x̂ = optimal_decision(farmer_problem)
3-element Array{Float64,1}:
 170.0
  80.0
 250.0
# Inspect optimal recourse decision in scenario 1
julia> optimal_recourse_decision(farmer_problem, 1)
6-element Array{Float64,1}:
    0.0
    0.0
  310.00000000000017
   48.000000000000036
 6000.0
    0.0
# Inspect optimal value
julia> objective_value(farmer_problem)
-108390.0
# Calculate expected value of perfect information
julia> EVPI(farmer_problem)
7015.6
# Calculate value of the stochastic solution
julia> VSS(farmer_problem)
1150.0
# Initialize with vertical structure
julia> farmer_ls = instantiate(farmer_model, [ξ₁,ξ₂,ξ₃], optimizer = LShaped.Optimizer);
# Set GLPK optimizer for the solving master problem
julia> set_optimizer_attribute(farmer_ls, MasterOptimizer(), GLPK.Optimizer);
# Set GLPK optimizer for the solving subproblems
julia> set_optimizer_attribute(farmer_ls, SubproblemOptimizer(), GLPK.Optimizer);
# Solve using L-shaped
julia> optimize!(farmer_ls)
L-Shaped Gap  Time: 0:00:00 (6 iterations)
  Objective:       -108390.0
  Gap:             0.0
  Number of cuts:  14
  Iterations:      6
# Initialize with horizontal structure
julia> farmer_ph = instantiate(farmer_model, [ξ₁,ξ₂,ξ₃],
                      optimizer = ProgressiveHedging.Optimizer);
# Set Ipopt optimizer for soving emerging subproblems
julia> set_optimizer_attribute(farmer_ph, SubproblemOptimizer(), Ipopt.Optimizer)
# Silence Ipopt
julia> set_optimizer_attribute(farmer_ph, RawSubproblemOptimizerParameter("print_level"), 0
# Solve using progressive-hedging
julia> optimize!(farmer_ph)
Progressive Hedging Time: 0:00:05 (86 iterations)
  Objective:   -108390.3601369591
  Primal gap:  3.984637579811031e-6
  Dual gap:    5.634811373041405e-6
  Iterations:  85
\end{lstlisting}

\subsection{Advanced solver configurations in the SPjl framework}
\label{sec:algor-impr}

The SPjl framework includes a variety of customizable improvements to the L-shaped and progressive-hedging algorithms. The possible variations of the classical algorithms included in the framework range from efficient implementations of influential research papers~\cite{ruszczynski_regularized_1986, linderoth_decomposition_2003, Fabian2006} to novel variants developed by the framework authors~\cite{cutaggregation} or others~\cite{adaptive_penalty, Wolf2013, Trukhanov2010}. We provide a summary of the improvements available for both L-shaped and progressive-hedging. In brief, each algorithm has a set of options that can be varied through a simple interface. In all examples, it is assumed that a given stochastic program instance \jlinl{sp} has been instantiated with an appropriate optimizer. We can then use \jlinl{set_optimizer_attribute(sp, option, value)} to customize the optimizer algorithm used by \jlinl{sp}.

\subsubsection{L-shaped}
\label{sec:l-shaped}

The L-shaped solver suite of SPjl includes a large set of customizable options. These are summarized below.

\textbf{Regularization}: A Regularization procedure limits the candidate search to a neighborhood of the current best iterate in the master problem. It tends to result in more effective cutting planes and improved performance of the L-shaped algorithm. Moreover, regularization enables warm-starting the L-shaped procedure with initial decisions. We have previously covered the regularization procedures in SPjl more in depth in~\cite{distlshaped}.

The SPjl framework includes the following regularizations: Trust-region regularization~\cite{linderoth_decomposition_2003}, Regularized decomposition~\cite{ruszczynski_regularized_1986}, Level set regularization~\cite{Fabian2006}. Since the two latter techniques involve solving problems with quadratic penalty terms, the SPjl framework also provide an option for replacing quadratic penalties with various linear approximations, if only a linear solver is available.
\sloppy

\textbf{Aggregation}: Cut aggregation can reduce communication overhead and load imbalance and yield major performance improvements in distributed settings. In the classical L-shaped algorithm~\cite{van_slyke_l-shaped_1969}, all cuts are aggregated every iteration. The authors of~\cite{Birge1988} suggested a multi-cut variant where cuts are added separately in a disaggregate form, which on average yields faster convergence. We recently explored a novel set of aggregation approaches~\cite{cutaggregation}, which are all included in SPjl.

\sloppy
\textbf{Consolidation}: Cut consolidation, as proposed by~\cite{Wolf2013}, is also implemented in SPjl to reduce load imbalance by removing stale cuts from the master.

\textbf{Execution}: In a distributed environment with multiple Julia processes, the execution policy of the L-shaped algorithm can be executed in a serial, synchronous or asynchronous mode. The synchronous variant runs the L-shaped algorithm in parallel using a map-reduce pattern each iteration. The asynchronous scheme is appropriate in a heterogeneous environment where some workers may finish slower than others. We show how these algorithm policies can be applied to increase performance on large-scale problems in Section~\ref{sec:numerical-benchmarks}.

\subsubsection{Progressive-hedging}
\label{sec:progressive-hedging}

\sloppy
The progressive-hedging solver suite shares a few options with the L-shaped suite. First, as each subproblem in the progressive-hedging procedure includes a quadratic penalty term, the same linear approximations as for L-shaped regularizations can be applied. Second, just like the L-shaped solvers, the progressive-hedging algorithms can be run serially, synchronously or asynchronously.

\textbf{Penalization}: The convergence rate of the progressive-hedging algorithm is sensitive to the choice of the penalty parameter $r$. The SPjl framework supports both a fixed penalty parameter and the adaptive strategy introduced in~\cite{adaptive_penalty}.

\section{Implementation details}
\label{sec:implementation}

In this section, we provide a summary of the main software innovations in SPjl. We also discuss the implementation of the framework's distributed capabilities. The inner workings of SPjl are primarily based on two ideas: deferred model instantiation and data injection. In brief, a model definition in SPjl is a recipe for how to use data structures when building optimization models, while the actual model creation is deferred until data is provided. When a specific model is instantiated, the provided data is injected where required to construct the model. The main effect of this approach is that the stochastic model formulation is separated from the design of stochastic data parameters, which makes the SPjl framework versatile and flexible to use. For instance, it is possible to test small instances of a model locally to ensure that it is properly defined, and then run the same model in a distributed environment with a large set of scenarios. Deferred model instances and data injection also play a large role when distributing stochastic program instances in memory.

\subsection{Deferred model instantiation}
\label{sec:deferr-model-instant}

\sloppy
The advantages of deferred model instantiation is a smaller memory footprint and the ability to create various structures that use the first- and second-stage recipes as building blocks in a clever way. Examples include the deterministic, vertical, and horizontal structures, as well as wait-and-see problems and expected-value problems. The technique is also a premise for implementing data injection. In contrast to standard JuMP models, SPjl models defined through the \jlinl{@stage} macros are not necessarily instantiated immediately. Instead, the user-defined Julia code that constructs the optimization problems is stored in lambda functions as model recipes. In other words, instead of creating and storing a JuMP object, the lines of code required to create the JuMP object is stored. This is achievable since Julia code is itself a data structure defined within the Julia language.

\sloppy
Deferred model instantiation is made possible through metaprogamming and the automatic reformulation bridges introduced in \jlinl{MathOptInterface}~\cite{moi}. These techniques allow us to add linking constraints between the stages that adhere to the data dependencies defined by the user. During model creation, any \jlinl{@decision} line in a \jlinl{@stage} definition creates special JuMP variables whose behaviour depends on the context of the instantiation. Any variable defined in this way can be included in \jlinl{@constraint} definitions in subsequent stages. See for example Listing~\ref{lst:simpledef}, where the last two constraint definitions in the second stage include references to \jlinl{x₁} and \jlinl{x₂} which were defined with \jlinl{@decision} in the first stage. Next, we will discuss in more detail how instantiation is implemented for the main underlying structures: deterministic, vertical, and horizontal. In addition, we explain how decision evaluation is implemented in the different structures.

\subsubsection{Deterministic structure}
\label{sec:determ-struct}

\sloppy
We construct the extensive form of a finite model~\eqref{eq:finitesp} in steps using the stored model recipes. First, we generate the first-stage model in full using the corresponding recipe. Next, we process all available scenarios iteratively. For each scenario, we apply the second-stage recipe and append the resulting subproblem to the extensive model. In this context, any variables defined with \jlinl{@decision} in the first stage are treated as regular JuMP variables. Before generating the subsequent scenario problem, we internally annotate the variables and constraints to associate them with the scenario they originated from. This labeling is visible in the printout shown in Listing~\ref{lst:simplediscrete}. During decision evaluation, all variables defined with \jlinl{@decision} are fixed to their corresponding values. The deterministic equivalent problem is then solved as usual, giving exactly~\eqref{eq:expectedresult}.

\subsubsection{Vertical structure}
\label{sec:vertical-struct}

\sloppy
The vertical structure, introduced in Section~\ref{sec:l-shaped-algorithm}, is also instantiated in steps. First, the first-stage master problem~\eqref{eq:lsmaster} is created using the corresponding recipe. Here, the \jlinl{@decision} variables are again treated as regular JuMP variables. Next, subproblem instances of the form~\eqref{eq:lssubprob} are created for each possible scenario using the second-stage recipe. During second-stage generation, first-stage variables annotated with \jlinl{@decision} enter the model as so called \emph{known decisions}. These are not optimization variables, but rather parameters with given values. This design reflects the fact that the first-stage decisions have already been taken when the second stage is reached. The values of the first-stage decisions can be entered into the second-stage constraints in which they appear through automatic reformulation bridges. Internally, all decisions defined in the first-stage are made known to the second stage by the \jlinl{@stochastic_model} macro. It is also possible to explicitly add \jlinl{@known} annotations to the second-stage definition to mark variables that originate from previous stages.

\sloppy
The subproblems are either stored in vector format on the master node or distributed on remote nodes as described in Section~\ref{sec:distr-comp}. We distribute new scenarios and generated subproblems as evenly as possible on remote nodes to achieve load balance. During decision evaluation, all variables defined with \jlinl{@decision} are fixed to their corresponding values in the first stage. Further, these values are communicated to all subproblems, that can then update their respective second-stage constraints. The first stage and second stage problems are then solved separately, in parallel if possible, and the results are map-reduced to form~\eqref{eq:expectedresult}.

\subsubsection{Horizontal structure}
\label{sec:horizontal-struct}

\sloppy
Instantiation of the horizontal structure introduced in Section~\ref{sec:progr-hedg-algor} is similar to instantiation of the vertical structure. They differ in that there is no master problem and in that the subproblems have the structure given in~\eqref{eq:phform} instead of~\eqref{eq:lssubprob}. The process for generating subproblems of this wait-and-see form is equivalent to one iteration of the finite extensive form generation. In short, the first-stage recipe is applied, followed by applying the second-stage recipe on the scenario data corresponding to the subproblem. Now, the variables defined with \jlinl{@decision} are again treated as standard JuMP variables. Note that generation of the expected-value-problem~\eqref{eq:ev} is equivalent to generating a wait-and-see model on the expected scenario of all available scenarios. The implementable solution $\xi$ that enter the horizontal form~\eqref{eq:phform} through the non-anticipative constraints is added as a known decision to the subproblems. During the progressive-hedging procedure, the value of $\xi$ can then be updated efficiently through bridges. This design is also used to implement proximal terms in the regularized variants of the L-shaped algorithm. Decision evaluation is performed similar to the other structures. In each subproblem, the first-stage decisions are fixed to their corresponding values and the subproblem is solved as usual. The results are then map-reduced to form~\eqref{eq:expectedresult}. Again, the decision evaluation process is embarrassingly parallel in a distributed environment.

\subsection{Data injection}
\label{sec:data-injection}

\sloppy
Data injection is the second software pattern used to separate model and data design in SPjl. The aim is to make an object independent of how its dependencies are created. In SPjl, the dependencies consist of the data required to construct the optimization problems as described by the model recipes. The data includes uncertain parameters, as well as first-stage decisions and deterministic parameters. By adopting this approach, users of SPjl can focus on the design of the optimization model and the uncertainty model separately, while the framework is responsible for combining these designs into actual stochastic program instances. In the following, we describe the data injection functionality in more detail.

When an SPjl model is formulated using \jlinl{@stochastic_model}, special annotations are used inside the \jlinl{@stage} blocks to specify points of data injection. These annotations inform the framework which parameters are necessary to construct the model according to the \jlinl{@stochastic_model} definition. The \jlinl{@stage} macro transforms the stage blocks into anonymous lambda functions that map supplied data into optimization problems. Internally, when the user wants to instantiate the defined SPjl model, the required data is passed to the stored lambda functions according to one of the instantiation procedures outlined in the previous section.

We give a short review of the different types of data dependencies that can be specified in an SPjl model. Consider the simple second-stage formulation in Listing~\ref{lst:injectionexample}, which includes several data injection annotations.
\begin{lstlisting}[language = julia, float, floatplacement=H, caption={Simple showcase of data injection in SPjl}, label={lst:injectionexample}]
@stochastic_model begin
    @stage 1 begin
        @decision(model, x)
    end
    @stage 2 begin
        @parameters d
        @known x
        @uncertain ξ
        @recourse(model, y <= d)
        @constraint(model, x + y <= ξ)
    end
end
\end{lstlisting}
\noindent
\sloppy

\textbf{Deterministic data}: The \jlinl{@parameters} annotation specifies scenario-independent data, i.e., deterministic parameters that are the same across all scenarios. Default parameter values can be specified inside the \jlinl{@parameters} block. Otherwise, the values must be supplied during instantiation.

\textbf{Uncertain data}: The \jlinl{@uncertain} annotation specifies the scenario-specific data. The scenario-dependent values are either created and supplied directly by the user or by some user-defined sampler object that models the uncertainty.

\textbf{Decisions}: The \jlinl{@known} annotation makes the first-stage decision \jlinl{x} available in the second-stage. Note again that \jlinl{@known} annotations are implicitly added by \jlinl{@stochastic_model} because of the \jlinl{@decision x} in the first stage. When the second-stage generator is run, the framework will have already created a decision variable \jlinl{x} using the first-stage generator, either as a standard JuMP variable or as a fixed known decision. All such first-stage variables are injected into the second-stage generator. These can then be used as if they were ordinary JuMP variables. See for example the last \jlinl{@constraint} definition in the second stage of Listing~\ref{lst:injectionexample}.

\textbf{Models}: The \jlinl{model} keyword is a placeholder for a JuMP object that stores the actual optimization problem. In a deterministic structure, the model object is the same in every generator call. In the block-decomposition structures, the generators are instead applied to multiple JuMP models that form subproblems.

\sloppy
The use of data injection adds versatility to the framework. The user is only restricted to use the \jlinl{model} keyword in the JuMP macro calls. Otherwise, all JuMP features are supported in the stage blocks. Also, there is no restriction on the scenario data types. Hence, instead of a simple structure with fields, it is possible to define a more complex data type that for example performs calculations at runtime to determine optimization parameters. In addition, any Julia methods defined on the scenario type become available in the stage blocks. This allows the user to design complex models of the uncertainty orthogonally to the definition of the stochastic program.

Because the model definition is decoupled from the data, it is possible to send the model recipe to a remote process where the scenario data resides and create the model from there. This is the foundation of the distributed implementation described next.

\subsection{Distributed computations}
\label{sec:distr-comp}

\sloppy
SPjl has distributed capabilities for both modeling, analysis and optimization. All implementations rely on the \jlinl{Distributed} module in Julia. This allows us to develop SPjl using high-level abstractions that utilize the efficient low-level communication protocols in Julia. In this way, the same codebase can be used to distribute computations locally, using shared-memory, and remotely, in a cluster or in the cloud.

\sloppy
Distributed computing in Julia is centered around the concepts of remote references and remote calls. Remote references are used to administer which node particular data resides on and to provide the remaining processes access to the remote data. Remote calls are used to schedule tasks on the nodes. Any process can \jlinl{wait} on a remote reference, which blocks until data can be fetched, and then \jlinl{fetch} the result when it is ready. The \jlinl{RemoteChannel} objects are special remote references where processes can also \jlinl{put!} data. Besides, specialized channel objects can be designed for specific data types. This feature is used frequently in the implementation of the distributed structured solvers.

\subsubsection{Distributed stochastic programs}
\label{sec:distr-stoch-progr}

\sloppy
The distributed capabilities of SPjl were designed with the aim to minimize the amount of data passing. This is mainly achieved through the deferred instantiation and data injection techniques outlined above. In principle, a stochastic program can be instantiated in a distributed environment by passing all necessary data to each worker node. However, the data injection technique is independent of the way data is created. Therefore, a far more efficient approach is to let the workers generate the necessary scenario data and the optimization models themselves, with minimal data passing. This is possible since SPjl has support for passing lightweight sampler objects capable of randomly generating scenario data, such as the one in defined in Listing~\ref{lst:simplecont}, along with passing the lightweight model recipes created in the \jlinl{@stage} blocks. Scenario data and subproblems can then be generated in parallel on the worker nodes. The master keeps track of the scenario distribution and ensures that new scenarios and subproblems are generated on available workers in a way that promotes load-balance.

If multiple Julia processes are available, then any instantiated stochastic program in SPjl is automatically distributed in memory according to either a vertical structure or a horizontal structure. In a vertical structure, the master node administers the first-stage problem and schedules tasks and data transfers. In a horizontal structure, the master node is only responsible for task scheduling and data transfers. Aside from distributing the models in memory, SPjl parallelizes as many computations as possible. In many cases, speedups stem from subtasks being embarrassingly parallel over the independent subproblems. For example, this occurs during decision evaluation and calculation of EVPI and VSS. In these instances, the master schedules the same computation tasks on all workers using remote calls and then initiates any necessary reductions after the workers have finished using a standard map-reduce pattern. The more involved parallelization strategies in SPjl relate mostly to the structure-exploiting distributed solvers described in more detail in Appendix~\ref{sec:distr-comp-deta}.

\subsubsection{Distributed structured optimization algorithms}
\label{sec:distr-comp-deta}

The implementations of the distributed structured solvers are also centered around remote calls and channels. Here, remote calls are used to initiate running tasks on every worker node, and the algorithm logic is driven by having the master and worker tasks wait on and write/fetch to/from specialized queue channels.

\sloppy
In the case of the L-shaped method, whenever the master node re-solves the master problem~\eqref{eq:lsmaster}, it writes the new decision vector to a specialized \jlinl{Decision} channel. It then sends a corresponding index to a \jlinl{Work} channel on every remote node. Every worker continuously fetches tasks from its \jlinl{Work} channel and uses the acquired index to fetch the latest decision vector from the master. Every new decision candidate infers a batch of subproblems to solve for each worker. After a worker has solved a subproblem~\eqref{eq:lssubprob}, it sends the computed cutting planes to a \jlinl{CutQueue} channel on the master. The master continuously fetches cuts from the \jlinl{CutQueue} and appends them to the master problem. In the synchronous variant, the master only updates after all workers have finished their work for the current iteration. In other words, the synchronous algorithm is driven by the master node initiating and waiting for worker tasks through remote calls. In the asynchronous version, the master updates after it has received $\kappa n$ cuts, where $n$ is the total number of subproblems. Timestamps are communicated throughout to keep track of the algorithm history and allow synchronized convergence checks. All subproblems are solved to completion each iteration regardless of the value of $\kappa$, to be able to check convergence properly. When the master has received all cuts corresponding to a specific iteration, it performs a convergence check and terminates if appropriate. For clarity, the procedure is illustrated in Fig.~\ref{fig:async-l-shaped}. A similar design is used to implement synchronous and asynchronous variants of the progressive-hedging algorithm.

\begin{figure}
  \centering
  \input{async-l-shaped}
  \caption{Asynchronous L-shaped procedure}
  \label{fig:async-l-shaped}
\end{figure}

\section{Numerical benchmarks}
\label{sec:numerical-benchmarks}

We now evaluate the distributed performance of SPjl by benchmarking the structure-exploiting solvers on a large-scale planning problem. The numerical experiments are performed in a multi-node setup where a laptop computer acts as the master node and a desktop compute server of up to $\num{32}$ cores provides worker nodes.

\subsection{The SSN problem}
\label{sec:ssn}

\begin{figure}
  \centering
  \input{ssn_confidence_intervals.tex}
  \caption{$90\%$ Confidence intervals around the optimal value of the SSN problem as a function of sample size.}
  \label{fig:ssn_confidence_intervals}
\end{figure}
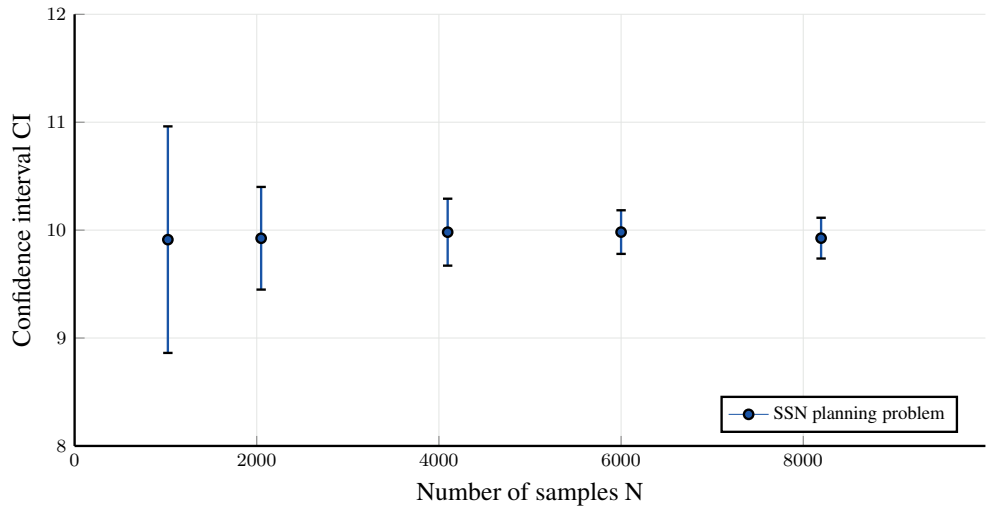

We evaluate the solvers on the telecommunications problem SSN, first introduced in~\cite{ssn}. This problem is often included in similar benchmarks~\cite{linderoth_decomposition_2003, Trukhanov2010}. The SSN problem is formulated to plan bandwidth capacity expansion in a network before customer demands are known. The problem is freely available in the SMPS format~\footnote{https://core.isrd.isi.edu/}. The problem has $\num{89}$ decision variables in the first stage, and $\num{706}$ variables and $\num{175}$ constraints in the second stage. We first run an SAA procedure to gauge the number of scenarios required to obtain a stable solution. The results are shown in Figure~\ref{fig:ssn_confidence_intervals}. There is no visible improvement after $\num{6000}$ scenarios. Moreover, the confidence interval around the optimal value is considered relatively tight at this point and is consistent with similar experiments~\cite{saacomp}. With $\num{6000}$ scenarios, the extensive form of the SAA model has 4.2 million variables and 5.3 million constraints, and about 20 minutes is required to build and solve the extensive form using Gurobi~\cite{gurobi}. From this baseline, we run the distributed benchmarks.

\subsubsection{Benchmarks}
\label{sec:distr-benchm}

We evaluate the structured solvers by solving distributed SSN instances of $\num{6000}$ scenarios. Benchmarks are performed using the Julia package BenchmarkTools.jl, which schedules multiple solve procedures and reports median computation times. Every solver runs until convergence criteria are reached with a relative tolerance of $10^{-2}$. The master node is a laptop computer with a 2.6 GHz Intel Core i7 processor and 16 GB of RAM. We spawn workers on a remote multi-core machine with two 3.1 GHz Intel Xeon processors (total 32 cores) and 128 GB of RAM. The two machines were $\num{30}$ kilometers apart at the time of the experiments. The time required to pass a single decision or optimality cut at this distance is about $\num{0.01}$ seconds. Hence, the communication latency is small, but not negligible as will be apparent in the results. For single-core experiments we only run the procedures once because the time to convergence is long and the measurement variance becomes relatively small. Throughout, the Gurobi optimizer~\cite{gurobi} is used to solve emerging subproblems.

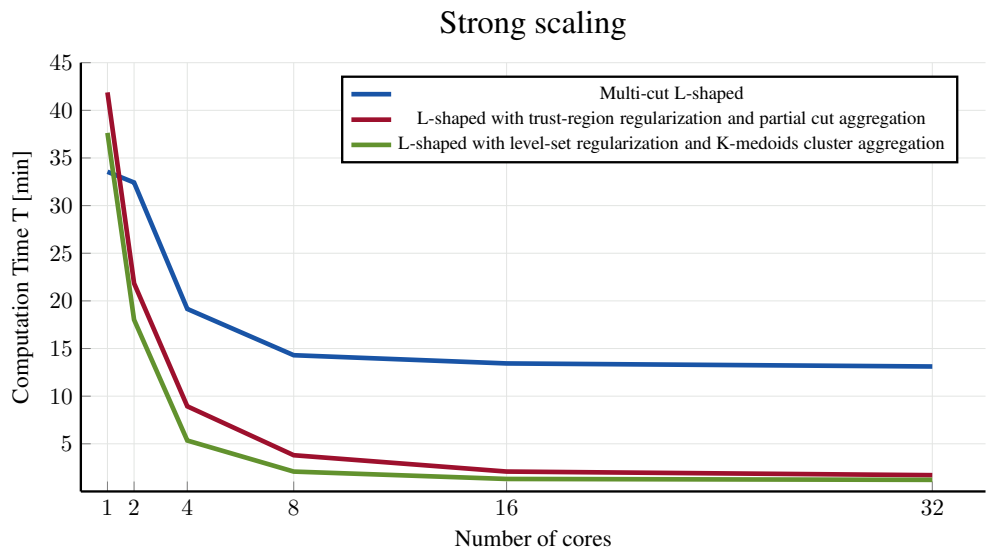
\begin{figure}
  \centering
  \input{strong_scaling}
  \caption{Median computation time required for different L-shaped algorithms to solve SSN instances of 6000 scenarios, as a function of number of worker cores. All experiments were run under synchronous execution.}
  \label{fig:strong_scaling}
\end{figure}

We first benchmark a set of L-shaped solvers. The nominal method is the multi-cut L-shaped algorithm without any advanced configuration. On average, this algorithm requires $\num{19}$ iterations and $\num{92000}$ optimality cuts to solve an SSN instance of $\num{6000}$ scenarios. This takes just over $\num{30}$ minutes on the master node under serial execution. We run a strong scaling test where the number of worker cores on the remote machine is doubled in size up to $\num{32}$ cores. Apart from multi-cut L-shaped, we also evaluate two variants with advanced algorithm policies. Specifically, one solver is configured to use trust-region regularization and partial cut aggregation with $\num{32}$ cuts in each bundle. This aggregation scheme is static; the cuts are partitioned into groups of $\num{32}$ in the same order each iteration. The second solver is configured to use level-set regularization and K-medoids cluster aggregation. This is a dynamic aggregation scheme where the cuts are clustered using the K-medoids algorithm based on a generalized cut distance matrix each iteration. We fix the partitioning scheme of the dynamic method after the first five iterations, as outline in~\cite{cutaggregation}. All solvers are configured to use synchronous execution. The results from these experiments are shown in Figure~\ref{fig:strong_scaling}.

First, we just consider the multi-cut method. The initial scaling is very poor with almost no speedup. We then observe speedups up to eight cores upon which the scaling curve flattens. The primary sources of inefficiency in distributed L-shaped algorithms are communication latency and load imbalance. This is especially true for multi-cut L-shaped because all cuts are passed separately and the master increases in size by the maximum number of constraints possible each iteration. We re-ran the two-core experiment on the master node with local threads as workers. In other words, without communication overhead. The time to convergence was then about $\num{18}$ minutes. With $\num{0.01}$ seconds required to pass a single cut and $\num{92000}$ cuts passed in total, this accounts for the extra $\num{15}$ minutes required to converge in the multi-node setup. Therefore, we can conclude that much of the inefficiency stems from communication latency. The fact that the scaling curve flattens stems mostly from load imbalance. In the final iterations, most of the time is spent solving the now large master problem or passing cuts, so the worker nodes are not utilized optimally.

Next, we consider the advanced methods. The distributed performance is significantly improved compared to the multi-cut method. The main reason for this is that cut aggregation reduces both communication latency and load imbalance. Because cuts are aggregated, less data is passed each iteration. Further, the master problem does not grow as fast. Hence, the workload is more evenly spread out between master and workers, which improves parallel performance. In this particular case, the more advanced aggregation scheme yields slightly better performance, but it could also hold that level-set regularization is more performant than trust-region regularization on the SSN problem. Even with cut aggregation, the size of the master eventually exceeds the size of the subproblems and data passing still becomes a bottle-neck as the number of cores increase. Therefore, the scaling curves still flatten for larger numbers of cores. We do not claim that these configurations are the best possible. We can for example note that they are not optimal for single-core execution where both variants are outperformed by the multi-cut method. Also, the parallel efficiency increase as workers are added is not uniform. This is because the aggregation schemes are more optimal for some work granularities. We could possibly improve the convergence times further by parameter tuning. For this particular configurations, we could also let a processor on the remote machine act as the master node and remove communication latency all together. However, we believe that our results are a strong encouragement for the distributed capability of the SPjl framework. With non-negligible communication latency we are able to solve a large-scale planning problem in just over a minute by employing some of the readily available algorithm policies in the framework. This can be seen as a proof of concept for running industrial planning problems in a modern cloud architecture.

We tested the algorithms with asynchronous execution as well, but saw no performance improvements. Even though there is communication latency between the master node and the remote node, worker performance is even. Moreover, the subproblems are equally difficult to solve. There is therefore no immediate gain from introducing asynchrony and the overhead from doing so decreases performance. The asynchronous variants are expected to yield better performance in a more heterogeneous environment with stalling workers.

\begin{figure}
  \centering
  \input{ph_scaling}
  \caption{Median computation time required for the progressive-hedging algorithm to solve SSN instances of 6000 scenarios, as a function of number of worker cores. All experiments were run under synchronous execution.}
  \label{fig:ph_scaling}
\end{figure}
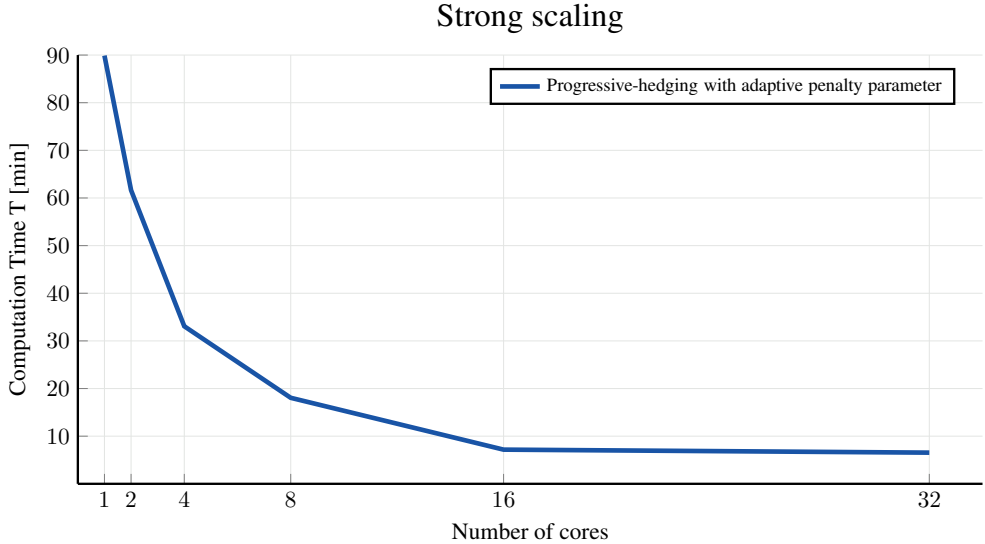

Next, we evaluate the performance of the progressive-hedging methods. Using the nominal method, we did not observe convergence even after long waiting times. Using the adaptive penalty policy eventually yields convergence. We configure the solvers to use adaptive penalty and synchronous execution and run the same strong scaling experiment as for the L-shaped methods. The results are shown in Figure~\ref{fig:ph_scaling}.

Although at much worse time-to-solution than the L-shaped methods initially, the distributed progressive-hedging algorithm displays great scaling and outperforms the multi-cut L-shaped method after $\num{16}$ cores. The efficiency probably stems from the problem being load-balanced across the workers. Communication latency again becomes a bottle-neck at $\num{32}$ cores from which we attribute the worsened scaling. Again, the subproblems appear equally difficult as there were no stalling workers. Consequently, we did not observe any speedups from running the asynchronous variant. The time to convergence is notably large and the progressive-hedging method is consistently outperformed by the advanced L-shaped methods. This is not surprising as we have spent more time on L-shaped improvements. Future work includes further algorithmic improvements to the progressive-hedging algorithms.

\section{Conclusion}
\label{sec:conclusion}

In this work, we have presented an open-source framework, \jlinl{StochasticPrograms.jl}, for large-scale stochastic programming. It is written entirely in Julia and includes both modeling tools and solver algorithms. The framework is designed for distributed computations and naturally scales to high-performance clusters or the cloud. By using the extensive form, which is efficiently generated using metaprogramming techniques, stochastic program instances can be solved using open-source or commercial solvers. Through deferred model instantiation, data injection, and clever algorithm policies, the framework can operate in distributed architectures with minimal data passing. In addition, several analysis tools and stochastic programming constructs are included with efficient implementations, many of which can run in parallel.

The framework also includes a solver suite of scalable algorithms that exploit the structure of the stochastic programs. The structured solvers are shown to perform well on large-scale planning problems. High parallel efficiency is achieved for distributed L-shaped methods using cut aggregation techniques and regularizations. Moreover, distributed progressive-hedging algorithms are accelerated using an adaptive penalty procedure. The solver suites are made modular through a policy-based design, so that future improvements can readily be added.

There are several directions for future additions to the framework. First, SPjl does not yet fully support multi-stage problems. We have finished an infrastructure for representing multi-stage problems in a way that leverages the two-stage design. Ongoing work involves designing a suitable Julian syntax for encoding transitive probabilities in a multi-stage scenario tree. Second, we will consider further algorithmic improvements to the existing L-shaped and progressive-hedging solvers. We also want to explore alternative sample-based approaches to the SAA method where the sampling is instead performed inside the structure-exploiting algorithm procedure. Examples of such approaches include L-shaped with importance sampling~\cite{Infanger1992} or stochastic decomposition~\cite{sd}.

\sloppy
The framework is well-tested through continuous integration and is freely available on Github\footnote{\url{https://github.com/martinbiel/StochasticPrograms.jl}}. A comprehensive documentation is included\footnote{\url{https://martinbiel.github.io/StochasticPrograms.jl/dev/}}. The modeling framework, \jlinl{StochasticPrograms.jl}, exists as a registered Julia package, which can be installed and run in any interactive Julia session.

\bibliographystyle{unsrt}
\bibliography{references}

\end{document}

%% file: async-l-shaped.tex
\resizebox{0.8\textwidth}{!}{
  \begin{tikzpicture}
    % Anchor
    \node at (0,2.5) (anchor1) {};

    % Decisions
    \node[scale=0.2,below of=anchor1,xshift=-100pt] (D) {$\mathcal{D}:$};
    \node[draw,scale=0.2, right of = D] (x0) {\color{kth-green}{$x_0$}};
    \node[draw,scale=0.2, right of = x0] (x1) {$\phantom{x_1}$};
    \node[draw,scale=0.2, right of = x1] (x2) {$\phantom{x_2}$};
    \node[scale=0.2, right of = x2] (x3) {$\cdots$};
    \node[draw,fit=(x0)(x1)(x2)(x3),inner sep=1pt] (Dbox) {};

    % Cuts
    \node[scale=0.2,below of = D] (C) {$\mathcal{C}:$};
    \node[draw,scale=0.2, right of = C] (c0) {$\phantom{c_0}$};
    \node[draw,scale=0.2, right of = c0] (c1) {$\phantom{c_1}$};
    \node[draw,scale=0.2, right of = c1] (c2) {$\phantom{c_2}$};
    \node[scale=0.2, right of = c2] (c3) {$\cdots$};
    \node[draw,fit=(c0)(c1)(c2)(c3),inner sep=1pt] (Cbox) {};

    % Master problem
    \node[scale=0.2,right of=Dbox,node distance=160pt] (masterprob) {$\begin{aligned}
        \minimize_{x \in \mathbb{R}^n} & \quad c^T x + \phantom{\sum_{s = 1}^{n} \theta_s} \\
        \text{s.t.} & \quad Ax = b \\
        & \quad \phantom{\partial \mathbf{Q} x + \theta_s \geq \mathbf{q}, \quad s = 1,\dots,n} \\
        & \quad x \geq 0
      \end{aligned}$};

    % Master
    \node[draw,fit=(D)(Dbox)(C)(Cbox)(masterprob),inner sep=1pt] (master) {};
    \node[above of=master,scale=0.4,node distance = 15pt] {\textbf{Master}};

    % Workers
    % Work 1
    \node[scale=0.2,above left=of anchor1] (work1) {$\mathcal{W}_1:$};
    \node[draw,scale=0.2, right of = work1] (w10) {$\color{kth-green}{1}$};
    \node[draw,scale=0.2, right of = w10] (w11) {$\phantom{2}$};
    \node[draw,scale=0.2, right of = w11] (w12) {$\phantom{3}$};
    \node[scale=0.2, right of = w12] (w13) {$\cdots$};
    \node[draw,fit=(w10)(w11)(w12)(w13),inner sep=1pt] (work1box) {};
    % Subproblems 1
    \node[scale=0.2,below of=work1] (S1) {$\mathcal{S}_1:$};
    \node[draw,scale=0.2,below of=work1box,node distance=50pt,xshift=13pt] (subprob1) {$\begin{aligned}
        \minimize_{y_i \in \mathrm{R}^m} & \quad q_s^T y_s \\
        \text{s.t.} & \quad Wy_s = h_s - T_s \color{kth-green}{x_0} \\
        & \quad y_s \geq 0
      \end{aligned}$};
    % W1
    \node[draw,fit=(work1)(work1box)(subprob1),inner sep=1pt] (W1) {};
    \node[above of=W1,scale=0.4,node distance = 15pt] {\textbf{Worker }$1$};
    % Work r
    \node[scale=0.2,above right=of anchor1] (workr) {$\mathcal{W}_r:$};
    \node[draw,scale=0.2, right of = workr] (wr0) {$\color{kth-green}{1}$};
    \node[draw,scale=0.2, right of = wr0] (wr1) {$\phantom{2}$};
    \node[draw,scale=0.2, right of = wr1] (wr2) {$\phantom{3}$};
    \node[scale=0.2, right of = wr2] (wr3) {$\cdots$};
    \node[draw,fit=(wr0)(wr1)(wr2)(wr3),inner sep=1pt] (workrbox) {};
    % Subproblems r
    \node[scale=0.2,below of=workr] (Sr) {$\mathcal{S}_r:$};
    \node[draw,scale=0.2,below of=workrbox,node distance=50pt,xshift=13pt] (subprobr) {$\begin{aligned}
        \minimize_{y_i \in \mathrm{R}^m} & \quad q_s^T y_s \\
        \text{s.t.} & \quad Wy_s = h_s - T_s \color{kth-green}{x_0} \\
        & \quad y_s \geq 0
      \end{aligned}$};
    % Wr
    \node[draw,fit=(workr)(workrbox)(subprobr),inner sep=1pt] (Wr) {};
    \node[above of=Wr,scale=0.4,node distance = 15pt] {\textbf{Worker }$r$};

    \node at ($(W1)!0.5!(Wr)$) {$\cdots$};

    % Send work
    \draw[kth-lightblue,->] (master) -- node[pos=0.3,scale=0.25,xshift=-20pt] {pass} (w10);
    \draw[kth-lightblue,->] (master) -- node[pos=0.3,scale=0.25,xshift=20pt] {pass} (wr0);

    % Get decision
    \draw[kth-pink,->] (x0) -- node[pos=0.5,scale=0.25,xshift=20pt] {fetch} (-0.2,3.25);
    \draw[kth-pink,->] (x0) -- node[pos=0.5,scale=0.25,xshift=25pt] {fetch} (2.15,3.25);

    \node [below of=master, node distance = 15pt, scale=0.3] {(a) Master sends task to workers. Workers fetch latest decision vector.};

    % Fig 2

    % Anchor
    \node at (0,0) (anchor2) {};

    % Decisions
    \node[scale=0.2,below of=anchor2,xshift=-100pt] (D) {$\mathcal{D}:$};
    \node[draw,scale=0.2, right of = D] (x0) {$x_0$};
    \node[draw,scale=0.2, right of = x0] (x1) {$\color{kth-green}{x_1}$};
    \node[draw,scale=0.2, right of = x1] (x2) {$\phantom{x_2}$};
    \node[scale=0.2, right of = x2] (x3) {$\cdots$};
    \node[draw,fit=(x0)(x1)(x2)(x3),inner sep=1pt] (Dbox) {};

    % Cuts
    \node[scale=0.2,below of = D] (C) {$\mathcal{C}:$};
    \node[draw,scale=0.2, right of = C] (c0) {$\phantom{c_0}$};
    \draw[line width=0.1mm] (-0.5,-0.44)--(-0.5,-0.36);
    \node[draw,scale=0.2, right of = c0] (c1) {$\phantom{c_1}$};
    \node[draw,scale=0.2, right of = c1] (c2) {$\phantom{c_2}$};
    \node[scale=0.2, right of = c2] (c3) {$\cdots$};
    \node[draw,fit=(c0)(c1)(c2)(c3),inner sep=1pt] (Cbox) {};

    % Master problem
    \node[scale=0.2,right of=Dbox,node distance=160pt] (masterprob) {$\begin{aligned}
        \minimize_{x \in \mathbb{R}^n} & \quad c^T x + \color{kth-green}{\sum_{s = 1}^{n} \theta_s} \\
        \text{s.t.} & \quad Ax = b \\
        & \quad \color{kth-green}{\partial \mathbf{Q} x + \theta_s \geq \mathbf{q}, \quad i = 1,\dots,n} \\
        & \quad x \geq 0
      \end{aligned}$};

    % Master
    \node[draw,fit=(D)(Dbox)(C)(Cbox)(masterprob),inner sep=1pt] (master) {};
    \node[above of=master,scale=0.4,node distance = 15pt] {\textbf{Master}};

    % Workers
    % Work 1
    \node[scale=0.2,above left=of anchor2] (work1) {$\mathcal{W}_1:$};
    \node[draw,scale=0.2, right of = work1] (w10) {$1$};
    \node[draw,scale=0.2, right of = w10] (w11) {$\color{kth-green}{2}$};
    \node[draw,scale=0.2, right of = w11] (w12) {$\phantom{3}$};
    \node[scale=0.2, right of = w12] (w13) {$\cdots$};
    \node[draw,fit=(w10)(w11)(w12)(w13),inner sep=1pt] (work1box) {};
    % Subproblems 1
    \node[scale=0.2,below of=work1] (S1) {$\mathcal{S}_1:$};
    \node[draw,scale=0.2,below of=work1box,node distance=50pt,xshift=13pt] (subprob1) {$\begin{aligned}
        \minimize_{y_1 \in \mathrm{R}^m} & \quad q_1^T y_1 \\
        \text{s.t.} & \quad Wy_1 = h_1 - T_1 x_0 \\
        & \quad y_1 \geq 0
      \end{aligned}$};
    % W1
    \node[draw,fit=(work1)(work1box)(subprob1),inner sep=1pt] (W1) {};
    \node[above of=W1,scale=0.4,node distance = 15pt] {\textbf{Worker }$1$};
    % Work r
    \node[scale=0.2,above right=of anchor2] (workr) {$\mathcal{W}_r:$};
    \node[draw,scale=0.2, right of = workr] (wr0) {$1$};
    \node[draw,scale=0.2, right of = wr0] (wr1) {$\color{kth-green}{2}$};
    \node[draw,scale=0.2, right of = wr1] (wr2) {$\phantom{3}$};
    \node[scale=0.2, right of = wr2] (wr3) {$\cdots$};
    \node[draw,fit=(wr0)(wr1)(wr2)(wr3),inner sep=1pt] (workrbox) {};
    % Subproblems r
    \node[scale=0.2,below of=workr] (Sr) {$\mathcal{S}_r:$};
    \node[draw,scale=0.2,below of=workrbox,node distance=50pt,xshift=13pt] (subprobr) {$\begin{aligned}
        \minimize_{y_r \in \mathrm{R}^m} & \quad q_r^T y_r \\
        \text{s.t.} & \quad Wy_r = h_r - T_r x_0 \\
        & \quad y_r \geq 0
      \end{aligned}$};
    % Wr
    \node[draw,fit=(workr)(workrbox)(subprobr),inner sep=1pt] (Wr) {};
    \node[above of=Wr,scale=0.4,node distance = 15pt] {\textbf{Worker }$r$};

    \node at ($(W1)!0.5!(Wr)$) {$\cdots$};

    % Send cut
    \draw[kth-lightblue,->] (W1) -- node[pos=0.3,scale=0.25,xshift=-20pt] {pass} (c0);
    \draw[kth-lightblue,->] (Wr) -- node[pos=0.15,scale=0.25,xshift=-15pt] {pass} (c0);

    % Send work
    \draw[kth-lightblue,->] (master) -- node[pos=0.3,scale=0.25,xshift=-20pt] {pass} (w11);
    \draw[kth-lightblue,->] (master) -- node[pos=0.3,scale=0.25,xshift=20pt] {pass} (wr1);

    % Iteration progress
    \draw[kth-green,fill] (-0.545,-0.43) rectangle (-0.51,-0.37);

    \node [below of=master, node distance = 15pt, scale=0.3] {(b) Workers solve subproblems and send cuts to master. Master problem re-solved after $\kappa n$ cuts have been received.};
    \node [below of=master, node distance = 19pt, scale=0.3] {Master sends new task to workers when a new decision vector is ready.};

    % Fig 3

    % Anchor
    \node at (0,-2.5) (anchor3) {};

    % Decisions
    \node[scale=0.2,below of=anchor3,xshift=-100pt] (D) {$\mathcal{D}:$};
    \node[draw,scale=0.2, right of = D] (x0) {$x_0$};
    \node[draw,scale=0.2, right of = x0] (x1) {$\color{kth-green}{x_1}$};
    \node[draw,scale=0.2, right of = x1] (x2) {$\phantom{x_2}$};
    \node[scale=0.2, right of = x2] (x3) {$\cdots$};
    \node[draw,fit=(x0)(x1)(x2)(x3),inner sep=1pt] (Dbox) {};

    % Cuts
    \node[scale=0.2,below of = D] (C) {$\mathcal{C}:$};
    \node[draw,scale=0.2, right of = C] (c0) {$\phantom{c_0}$};
    \node[draw,scale=0.2, right of = c0] (c1) {$\phantom{c_1}$};
    \node[draw,scale=0.2, right of = c1] (c2) {$\phantom{c_2}$};
    \node[scale=0.2, right of = c2] (c3) {$\cdots$};
    \node[draw,fit=(c0)(c1)(c2)(c3),inner sep=1pt] (Cbox) {};

    % Master problem
    \node[scale=0.2,right of=Dbox,node distance=160pt] (masterprob) {$\begin{aligned}
        \minimize_{x \in \mathbb{R}^n} & \quad c^T x + \sum_{s = 1}^{n} \theta_s \\
        \text{s.t.} & \quad Ax = b \\
        & \quad \partial \mathbf{Q} x + \theta_s \geq \mathbf{q}, \quad i = 1,\dots,n \\
        & \quad x \geq 0
      \end{aligned}$};

    % Master
    \node[draw,fit=(D)(Dbox)(C)(Cbox)(masterprob),inner sep=1pt] (master) {};
    \node[above of=master,scale=0.4,node distance = 15pt] {\textbf{Master}};

    % Workers
    % Work 1
    \node[scale=0.2,above left=of anchor3] (work1) {$\mathcal{W}_1:$};
    \node[draw,scale=0.2, right of = work1] (w10) {$1$};
    \node[draw,scale=0.2, right of = w10] (w11) {$2$};
    \node[draw,scale=0.2, right of = w11] (w12) {$\phantom{3}$};
    \node[scale=0.2, right of = w12] (w13) {$\cdots$};
    \node[draw,fit=(w10)(w11)(w12)(w13),inner sep=1pt] (work1box) {};
    % Subproblems 1
    \node[scale=0.2,below of=work1] (S1) {$\mathcal{S}_1:$};
    \node[draw,scale=0.2,below of=work1box,node distance=50pt,xshift=13pt] (subprob1) {$\begin{aligned}
        \minimize_{y_1 \in \mathrm{R}^m} & \quad q_1^T y_1 \\
        \text{s.t.} & \quad Wy_1 = h_1 - T_1 x_0 \\
        & \quad y_1 \geq 0
      \end{aligned}$};
    % W1
    \node[draw,fit=(work1)(work1box)(subprob1),inner sep=1pt] (W1) {};
    \node[above of=W1,scale=0.4,node distance = 15pt] {\textbf{Worker }$1$};
    % Work r
    \node[scale=0.2,above right=of anchor3] (workr) {$\mathcal{W}_r:$};
    \node[draw,scale=0.2, right of = workr] (wr0) {$1$};
    \node[draw,scale=0.2, right of = wr0] (wr1) {$2$};
    \node[draw,scale=0.2, right of = wr1] (wr2) {$\phantom{3}$};
    \node[scale=0.2, right of = wr2] (wr3) {$\cdots$};
    \node[draw,fit=(wr0)(wr1)(wr2)(wr3),inner sep=1pt] (workrbox) {};
    % Subproblems r
    \node[scale=0.2,below of=workr] (Sr) {$\mathcal{S}_r:$};
    \node[draw,scale=0.2,below of=workrbox,node distance=50pt,xshift=13pt] (subprobr) {$\begin{aligned}
        \minimize_{y_r \in \mathrm{R}^m} & \quad q_r^T y_r \\
        \text{s.t.} & \quad Wy_r = h_r - T_r \color{kth-green}{x_1} \\
        & \quad y_r \geq 0
      \end{aligned}$};
    % Wr
    \node[draw,fit=(workr)(workrbox)(subprobr),inner sep=1pt] (Wr) {};
    \node[above of=Wr,scale=0.4,node distance = 15pt] {\textbf{Worker }$r$};

    \node at ($(W1)!0.5!(Wr)$) {$\cdots$};

    % Send cut
    \draw[kth-lightblue,->] (W1) -- node[pos=0.3,scale=0.25,xshift=-20pt] {pass} (c0);

    % Get decision
    \draw[kth-pink,->] (x0) -- node[pos=0.5,scale=0.25,xshift=25pt] {fetch} (2.15,-1.75);

    % Iteration progress
    \draw[kth-lightgreen,fill] (-0.545,-2.93) rectangle (-0.46,-2.87);

    % Convergence
    \node[scale=0.2,xshift=-50pt,fill,kth-lightgreen] {\phantom{$|Q - \Theta| \leq \tau(\epsilon+|Q|)$?}};
    \node[scale=0.2,xshift=-50pt] {$|Q - \Theta| \leq \tau(\epsilon+|Q|)$?};

    \node [below of=master, node distance = 15pt, scale=0.3] {(c) Convergence check when all cuts have been received. Ready workers fetch latest decision. Procedure continues.};
  \end{tikzpicture}
}

%% file: ssn_confidence_intervals.tex
\resizebox{0.8\textwidth}{!}{
  \begin{tikzpicture}[]
    \begin{axis}[height = {80.6mm}, legend pos = {south east}, ylabel = {Confidence interval CI}, xmin = {0}, xmax = {10000}, ymax = {12}, xlabel = {Number of samples N}, unbounded coords=jump,scaled x ticks = false,xlabel style = {font = {\fontsize{11 pt}{14.3 pt}\selectfont}, rotate = 0.0},xmajorgrids = true,xtick = {0.0,2000.0,4000.0,6000.0,8000.0},xticklabels = {$0$,$2000$,$4000$,$6000$,$8000$},xtick align = inside,xticklabel style = {font = {\fontsize{8 pt}{10.4 pt}\selectfont}, rotate = 0.0},x grid style = {color = kth-lightgray,
        line width = 0.25,
        solid},axis x line* = left,x axis line style = {line width = 1,
        solid},scaled y ticks = false,ylabel style = {font = {\fontsize{11 pt}{14.3 pt}\selectfont}, rotate = 0.0},ymajorgrids = true,ytick = {8,9,10,11,12},yticklabels = {$8$,$9$,$10$,$11$,$12$},ytick align = inside,yticklabel style = {font = {\fontsize{8 pt}{10.4 pt}\selectfont}, rotate = 0.0},y grid style = {color = kth-lightgray,
        line width = 0.25,
        solid},axis y line* = left,y axis line style = {line width = 1,
        solid},    xshift = 0.0mm,
      yshift = 0.0mm,
      title style = {font = {\fontsize{14 pt}{18.2 pt}\selectfont}, rotate = 0.0},legend style = {line width = 1,
        solid,font = {\fontsize{8 pt}{10.4 pt}\selectfont}},colorbar style={title=}, ymin = {8}, width = {152.40000000000003mm}]\addplot+[draw=none, color = kth-blue,
      line width = 0,
      solid,mark = *,
      mark size = 2.0,
      mark options = {
        color = black,
        fill = kth-blue,
        line width = 1,
        rotate = 0,
        solid
      }] coordinates {
        (1024.0, 9.911650607192012)
        (2048.0, 9.924494134526144)
        (4096.0, 9.980931270623948)
        (6000.0, 9.981725190294481)
        (8196.0, 9.925619697185482)
      };
      \addlegendentry{SSN planning problem}
      \addplot+ [color = kth-blue,
      line width = 1,
      solid,mark = -,
      mark size = 2.0,
      mark options = {
        color = black,
        line width = 1,
        rotate = 0,
        solid
      },forget plot]coordinates {
        (1024.0, 8.861984732842611)
        (1024.0, 10.961316481541415)
      };
      \addplot+ [color = kth-blue,
      line width = 1,
      solid,mark = -,
      mark size = 2.0,
      mark options = {
        color = black,
        line width = 1,
        rotate = 0,
        solid
      },forget plot]coordinates {
        (2048.0, 9.448876125955223)
        (2048.0, 10.400112143097065)
      };
      \addplot+ [color = kth-blue,
      line width = 1,
      solid,mark = -,
      mark size = 2.0,
      mark options = {
        color = black,
        line width = 1,
        rotate = 0,
        solid
      },forget plot]coordinates {
        (4096.0, 9.6707800103331)
        (4096.0, 10.291082530914798)
      };
      \addplot+ [color = kth-blue,
      line width = 1,
      solid,mark = -,
      mark size = 2.0,
      mark options = {
        color = black,
        line width = 1,
        rotate = 0,
        solid
      },forget plot]coordinates {
        (6000.0, 9.779612761905298)
        (6000.0, 10.183837618683665)
      };
      \addplot+ [color = kth-blue,
      line width = 1,
      solid,mark = -,
      mark size = 2.0,
      mark options = {
        color = black,
        line width = 1,
        rotate = 0,
        solid
      },forget plot]coordinates {
        (8196.0, 9.736661207365765)
        (8196.0, 10.114578187005199)
      };
    \end{axis}

  \end{tikzpicture}

}

%% file: strong_scaling.tex
\resizebox{0.8\textwidth}{!}{
  \begin{tikzpicture}[]
    \begin{axis}[height = {80.6mm}, legend pos = {north east}, ylabel = {Computation Time T [min]}, title = {Strong scaling}, xmin = {0}, xmax = {34}, ymax = {45}, xlabel = {Number of cores}, unbounded coords=jump,scaled x ticks = false,xlabel style = {font = {\fontsize{10 pt}{13.0 pt}\selectfont}, rotate = 0.0},xmajorgrids = true,xtick = {1.0,2.0,4.0,8.0,16.0,32.0},xticklabels = {$1$,$2$,$4$,$8$,$16$,$32$},xtick align = inside,xticklabel style = {font = {\fontsize{10 pt}{13.0 pt}\selectfont}, rotate = 0.0},x grid style = {color = kth-lightgray,
        line width = 0.25,
        solid},axis lines* = left,x axis line style = {line width = 1,
        solid},scaled y ticks = false,ylabel style = {font = {\fontsize{10 pt}{13.0 pt}\selectfont}, rotate = 0.0},ymajorgrids = true,ytick = {5, 10, 15, 20, 25, 30, 35, 40, 45},yticklabels = {$5$, $10$, $15$, $20$, $25$, $30$, $35$, $40$, $45$},ytick align = inside,yticklabel style = {font = {\fontsize{10 pt}{13.0 pt}\selectfont}, rotate = 0.0},y grid style = {color = kth-lightgray,
        line width = 0.25,
        solid},axis lines* = left,y axis line style = {line width = 1,
        solid},    xshift = 0.0mm,
      yshift = 0.0mm,
      ,title style = {font = {\fontsize{14 pt}{18.2 pt}\selectfont}, rotate = 0.0}
      ,legend style = {line width = 1,
        solid,font = {\fontsize{8 pt}{10.4 pt}\selectfont}},colorbar style={title=}, ymin = {0}, width = {152.4mm}]
      \addplot+ [color = kth-blue,
      line width = 2,
      solid,mark = none,
      mark size = 2.0,
      mark options = {
        fill = kth-blue,
        line width = 1,
        rotate = 0,
        solid
      }]coordinates {
        (1.0, 33.55)
        (2.0, 32.42417315311666)
        (4.0, 19.156102492716666)
        (8.0, 14.297370235699999)
        (16.0, 13.4350466065)
        (32.0, 13.11380154975)
      };
      \addlegendentry{Multi-cut L-shaped}
      \addplot+ [color = kth-red,
      line width = 2,
      solid,mark = none,
      mark size = 2.0,
      mark options = {
        fill = kth-red,
        line width = 1,
        rotate = 0,
        solid
      }]coordinates {
        (1.0, 41.8833333333333)
        (2.0, 21.845425012433335)
        (4.0, 8.935931102900001)
        (8.0, 3.8018041197000003)
        (16.0, 2.0870440991166666)
        (32.0, 1.70713432055)
      };
      \addlegendentry{L-shaped with trust-region regularization and partial cut aggregation}
      \addplot+ [color = kth-green,
      line width = 2,
      solid,mark = none,
      mark size = 2.0,
      mark options = {
        fill = kth-green,
        line width = 1,
        rotate = 0,
        solid
      }]coordinates {
        (1.0, 37.63333333333333)
        (2.0, 18.022236241766667)
        (4.0, 5.344734589216667)
        (8.0, 2.08504669955)
        (16.0, 1.3104172107166667)
        (32.0, 1.2157116196333333)
      };
      \addlegendentry{L-shaped with level-set regularization and K-medoids cluster aggregation}
    \end{axis}
  \end{tikzpicture}
}

%% file: ph_scaling.tex
\resizebox{0.8\textwidth}{!}{
  \begin{tikzpicture}[]
    \begin{axis}[height = {80.6mm}, legend pos = {north east}, ylabel = {Computation Time T [min]}, title = {Strong scaling}, xmin = {0}, xmax = {34}, ymax = {90}, xlabel = {Number of cores}, unbounded coords=jump,scaled x ticks = false,xlabel style = {font = {\fontsize{10 pt}{13.0 pt}\selectfont}, rotate = 0.0},xmajorgrids = true,xtick = {1.0,2.0,4.0,8.0,16.0,32.0},xticklabels = {$1$,$2$,$4$,$8$,$16$,$32$},xtick align = inside,xticklabel style = {font = {\fontsize{10 pt}{13.0 pt}\selectfont}, rotate = 0.0},x grid style = {color = kth-lightgray,
        line width = 0.25,
        solid},axis lines* = left,x axis line style = {line width = 1,
        solid},scaled y ticks = false,ylabel style = {font = {\fontsize{10 pt}{13.0 pt}\selectfont}, rotate = 0.0},ymajorgrids = true,ytick = {10, 20, 30, 40, 50, 60, 70, 80, 90},yticklabels = {$10$, $20$, $30$, $40$, $50$, $60$, $70$, $80$, $90$},ytick align = inside,yticklabel style = {font = {\fontsize{10 pt}{13.0 pt}\selectfont}, rotate = 0.0},y grid style = {color = kth-lightgray,
        line width = 0.25,
        solid},axis lines* = left,y axis line style = {line width = 1,
        solid},    xshift = 0.0mm,
      yshift = 0.0mm,
      ,title style = {font = {\fontsize{14 pt}{18.2 pt}\selectfont}, rotate = 0.0}
      ,legend style = {line width = 1,
        solid,font = {\fontsize{8 pt}{10.4 pt}\selectfont}},colorbar style={title=}, ymin = {0}, width = {152.4mm}]
      \addplot+ [color = kth-blue,
      line width = 2,
      solid,mark = none,
      mark size = 2.0,
      mark options = {
        fill = kth-blue,
        line width = 1,
        rotate = 0,
        solid
      }]coordinates {
        (1.0, 89.88843420491666)
        (2.0, 61.634777444625)
        (4.0, 33.072543491375)
        (8.0, 18.0547831021)
        (16.0, 7.1735764085166664)
        (32.0, 6.543205914766667)
      };
      \addlegendentry{Progressive-hedging with adaptive penalty parameter}
    \end{axis}
  \end{tikzpicture}
}

%% file: paper.bbl
\begin{thebibliography}{10}

\bibitem{Birge2011}
John~R. Birge and Fran{\c{c}}ois Louveaux.
\newblock {\em Introduction to Stochastic Programming}.
\newblock Springer New York, 2011.

\bibitem{Fleten2007}
Stein-Erik Fleten and Trine~Krogh Kristoffersen.
\newblock Stochastic programming for optimizing bidding strategies of a nordic
  hydropower producer.
\newblock {\em European Journal of Operational Research}, 181(2):916--928,
  2007.

\bibitem{Groewe-Kuska2005}
Nicole Gröwe-Kuska and Werner Römisch.
\newblock Stochastic unit commitment in hydrothermal power production planning.
\newblock In {\em Applications of Stochastic Programming}, pages 633--653.
  Society for Industrial and Applied Mathematics, 2005.

\bibitem{petra_real-time_2014}
C.~G. Petra, O.~Schenk, and M.~Anitescu.
\newblock Real-{Time} {Stochastic} {Optimization} of {Complex} {Energy}
  {Systems} on {High}-{Performance} {Computers}.
\newblock {\em Computing in Science Engineering}, 16(5):32--42, 2014.

\bibitem{Krokhmal2005}
P.~Krokhmal, S.~Uryasev, and G.~Zrazhevsky.
\newblock Numerical comparison of conditional value-at-risk and conditional
  drawdown-at-risk approaches: Application to hedge funds.
\newblock In {\em Applications of Stochastic Programming}, pages 609--631.
  Society for Industrial and Applied Mathematics, 2005.

\bibitem{Zenios2005}
Stavros~A. Zenios.
\newblock Optimization models for structuring index funds.
\newblock In {\em Applications of Stochastic Programming}, pages 471--501.
  Society for Industrial and Applied Mathematics, 2005.

\bibitem{Powell1987}
Warren~B. Powell.
\newblock An operational planning model for the dynamic vehicle allocation
  problem with uncertain demands.
\newblock {\em Transportation Research Part B: Methodological}, 21(3):217--232,
  1987.

\bibitem{Powell2005}
Warren~B. Powell and Huseyin Topaloglu.
\newblock Fleet management.
\newblock In {\em Applications of Stochastic Programming}, pages 185--215.
  Society for Industrial and Applied Mathematics, 2005.

\bibitem{glpk}
Andrew Makhorin.
\newblock Gnu linear programming kit, 2020.
\newblock \url{https://www.gnu.org/software/glpk/}.

\bibitem{gurobi}
{Gurobi Optimization}.
\newblock Gurobi optimizer reference manual, 2020.
\newblock \url{http://www.gurobi.com}.

\bibitem{Rockafellar1991}
R.~T. Rockafellar and Roger J.-B. Wets.
\newblock Scenarios and policy aggregation in optimization under uncertainty.
\newblock {\em Mathematics of Operations Research}, 16(1):119--147, 1991.

\bibitem{van_slyke_l-shaped_1969}
R.~Van~Slyke and Roger J.-B Wets.
\newblock L-{Shaped} {Linear} {Programs} with {Applications} to {Optimal}
  {Control} and {Stochastic} {Programming}.
\newblock {\em SIAM Journal on Applied Mathematics}, 17(4):638--663, 1969.

\bibitem{Bezanson2017}
Jeff Bezanson, Alan Edelman, Stefan Karpinski, and Viral~B. Shah.
\newblock Julia: A fresh approach to numerical computing.
\newblock {\em {SIAM} Review}, 59(1):65--98, 2017.

\bibitem{polojl}
Martin Biel, Arda Aytekin, and Mikael Johansson.
\newblock {POLO}.jl: Policy-based optimization algorithms in {J}ulia.
\newblock {\em Advances in Engineering Software}, 136:102695, 2019.

\bibitem{Dunning2017}
Iain Dunning, Joey Huchette, and Miles Lubin.
\newblock {JuMP}: A modeling language for mathematical optimization.
\newblock {\em {SIAM} Review}, 59(2):295--320, 2017.

\bibitem{moi}
Benoit Legat, Oscar Dowson, Joaquim~Dias Garcia, and Miles Lubin.
\newblock Mathoptinterface: a data structure for mathematical optimization
  problems.
\newblock {\em arXiv preprint arXiv:2002.03447}, 2020.

\bibitem{Watson2012}
Jean-Paul Watson, David~L. Woodruff, and William~E. Hart.
\newblock {PySP}: modeling and solving stochastic programs in python.
\newblock {\em Mathematical Programming Computation}, 4(2):109--149, 2012.
\newblock Cited on p. 129.

\bibitem{Hart2017}
William~E. Hart, Carl~D. Laird, Jean-Paul Watson, David~L. Woodruff, Gabriel~A.
  Hackebeil, Bethany~L. Nicholson, and John~D. Siirola.
\newblock {\em Pyomo {\textemdash} Optimization Modeling in Python}.
\newblock Springer International Publishing, 2017.

\bibitem{fortsp}
Francis Ellison, Gautam Mitra, Chandra Poojari, and Victor Zverovich.
\newblock Fortsp: A stochastic programming solver.
\newblock \url{http://www.optirisk-systems.com/manuals/FortspManual.pdf}, 2009.

\bibitem{Huchette2014}
Joey Huchette, Miles Lubin, and Cosmin Petra.
\newblock Parallel algebraic modeling for stochastic optimization.
\newblock In {\em 2014 First Workshop for High Performance Technical Computing
  in Dynamic Languages}. {IEEE}, 2014.

\bibitem{lubin_parallel_2013}
Miles Lubin, J.~A.~Julian Hall, Cosmin~G. Petra, and Mihai Anitescu.
\newblock Parallel distributed-memory simplex for large-scale stochastic {LP}
  problems.
\newblock {\em Computational Optimization and Applications}, 55(3):571--596,
  2013.

\bibitem{saa}
Wai-Kei Mak, David~P. Morton, and R.Kevin Wood.
\newblock Monte carlo bounding techniques for determining solution quality in
  stochastic programs.
\newblock {\em Operations Research Letters}, 24(1):47 -- 56, 1999.

\bibitem{saaconvergence}
Alan~J. King and Roger J.-B. Wets.
\newblock Epi-consistency of convex stochastic programs.
\newblock {\em Stochastics and Stochastics Reports}, 34(1-2):83--92, 1991.

\bibitem{saadist}
Alexander Shapiro.
\newblock Asymptotic analysis of stochastic programs.
\newblock {\em Annals of Operations Research}, 30(1):169--186, 1991.

\bibitem{saacomp}
Jeff Linderoth, Alexander Shapiro, and Stephen Wright.
\newblock The empirical behavior of sampling methods for stochastic
  programming.
\newblock {\em Annals of Operations Research}, 142(1):215--241, 2006.

\bibitem{Birge1988}
John~R. Birge and Fran{\c{c}}ois~V. Louveaux.
\newblock A multicut algorithm for two-stage stochastic linear programs.
\newblock {\em European Journal of Operational Research}, 34(3):384--392, 1988.

\bibitem{Rockafellar1976}
R.~T. Rockafellar.
\newblock Monotone operators and the proximal point algorithm.
\newblock {\em {SIAM} Journal on Control and Optimization}, 14(5):877--898,
  1976.

\bibitem{ruszczynski_regularized_1986}
Andrzej Ruszczyński.
\newblock A regularized decomposition method for minimizing a sum of polyhedral
  functions.
\newblock {\em Mathematical Programming}, 35(3):309--333, 1986.

\bibitem{linderoth_decomposition_2003}
Jeff Linderoth and Stephen Wright.
\newblock Decomposition {Algorithms} for {Stochastic} {Programming} on a
  {Computational} {Grid}.
\newblock {\em Computational Optimization and Applications}, 24(2-3):207--250,
  2003.

\bibitem{Fabian2006}
Csaba~I. F{\'{a}}bi{\'{a}}n and Zolt{\'{a}}n Sz{\H{o}}ke.
\newblock Solving two-stage stochastic programming problems with level
  decomposition.
\newblock {\em Computational Management Science}, 4(4):313--353, 2006.

\bibitem{cutaggregation}
Martin Biel and Mikael Johansson.
\newblock Dynamic cut aggregation in {L}-shaped algorithms.
\newblock {\em arXiv preprint arXiv:1910.13752}, 2019.
\newblock Submitted for consideration to the European Journal of Operational
  Research. Under review.

\bibitem{adaptive_penalty}
Shohre Zehtabian and Fabian Bastin.
\newblock {\em Penalty parameter update strategies in progressive hedging
  algorithm}.
\newblock CIRRELT, 2016.

\bibitem{Wolf2013}
Christian Wolf and Achim Koberstein.
\newblock Dynamic sequencing and cut consolidation for the parallel hybrid-cut
  nested l-shaped method.
\newblock {\em European Journal of Operational Research}, 230(1):143--156,
  2013.

\bibitem{Trukhanov2010}
Svyatoslav Trukhanov, Lewis Ntaimo, and Andrew Schaefer.
\newblock Adaptive multicut aggregation for two-stage stochastic linear
  programs with recourse.
\newblock {\em European Journal of Operational Research}, 206(2):395--406,
  2010.

\bibitem{distlshaped}
Martin Biel and Mikael Johansson.
\newblock Distributed {L}-shaped algorithms in {J}ulia.
\newblock In {\em 2018 {IEEE}/{ACM} Parallel Applications Workshop,
  Alternatives To {MPI} ({PAW}-{ATM})}. {IEEE}, 2018.

\bibitem{ssn}
Suvrajeet Sen, Robert~D. Doverspike, and Steve Cosares.
\newblock Network planning with random demand.
\newblock {\em Telecommunication Systems}, 3(1):11--30, 1994.

\bibitem{Infanger1992}
Gerd Infanger.
\newblock Monte carlo (importance) sampling within a benders decomposition
  algorithm for stochastic linear programs.
\newblock {\em Annals of Operations Research}, 39(1):69--95, 1992.

\bibitem{sd}
Julia~L. Higle and Suvrajeet Sen.
\newblock Stochastic decomposition: An algorithm for two-stage linear programs
  with recourse.
\newblock {\em Mathematics of Operations Research}, 16(3):650--669, 1991.

\end{thebibliography}
